\documentclass[10pt]{amsart}
\usepackage{graphicx}
\usepackage{cite,amssymb,amsmath}
\usepackage{amsthm}

\theoremstyle{plain}

\newtheorem{theorem}{Theorem}[section]
\newtheorem{proposition}[theorem]{Proposition}
\newtheorem{corollary}[theorem]{Corollary}
\newtheorem{conjecture}[theorem]{Conjecture}
\newtheorem{lemma}[theorem]{Lemma}

\theoremstyle{remark}
\newtheorem*{remark}{Remark}

\newtheorem*{note}{Note}

\newcommand{\newoperator}[2]{\DeclareMathOperator{#1}{#2}}

\newoperator{\AH}{AH}
\newoperator{\AJ}{AJ}
\newoperator{\BS}{BS}               
\newoperator{\GF}{GF}               
\newoperator{\GH}{GH}
\newoperator{\GJ}{GJ}
\newoperator{\Isom}{Isom}
\newoperator{\KS}{KS}               
\newoperator{\ML}{ML}
\newoperator{\Nbhd}{Nbhd}
\newoperator{\PL}{PL}
\newoperator{\QF}{QF}
\newoperator{\QH}{QH}
\newoperator{\QJ}{QJ}
\newoperator{\QpH}{Q^\prime H}
\newoperator{\Rad}{Rad}
\newoperator{\Tr}{Tr}
\newoperator{\Vol}{Vol}
\newoperator{\area}{area}
\newoperator{\base}{base}
\newoperator{\codim}{codim}
\newoperator{\const}{const}
\newoperator{\degree}{degree}       
\newoperator{\glue}{glue}
\newoperator{\graft}{graft}
\newoperator{\inj}{inj}             
\newoperator{\interior}{int}
\newoperator{\length}{length}
\newoperator{\llength}{\underline length}
\newoperator{\ls}{ls}
\newoperator{\qf}{qf}
\newoperator{\skin}{skin}
\newoperator{\teich}{\textsl{T}}
\newoperator{\unglue}{unglue}
\newoperator{\wb}{wb}               
\newoperator{\window}{window}
\newoperator{\ws}{ws}               

\newcommand{\hy}{\mathbb{H}}
\newcommand{\integers}{\mathbb{Z}}
\newcommand{\complexes}{\mathbb{C}}
\newcommand{\proj}{\mathbb{P}}
\newcommand{\reals}{\mathbb{R}}
\newcommand{\euclidean}{\mathbb{E}}

\newcommand{\refto}[1]{Ref.~\citen{#1}}
\newcommand{\refin}[2]{#1 of Ref.~\citen{#2}}
\newcommand{\fund}[1]{\pi_1(#1)}
\newcommand{\set}[1]{\{#1\}}

\newcommand{\Chat}{\hat{\complexes}}

\newcommand{\ah}{\AH(M)}
\newcommand{\arrow}{\rightarrow}
\newcommand{\boundary}{\partial}

\newcommand{\cusps}{{\hbox{\rm cusps}}}

\newcommand{\homotopic}{\simeq}

\newcommand{\inverse}{^{-1}}

\newcommand{\sinfty}{S_\infty^2}
\newcommand{\thin}{{\hbox{\rm thin}}}
\newcommand{\tps}{\Chat-\set{0,1,\infty}} 

\newcommand{\zeemod}{\integers_}

\begin{document}
\title{Hyperbolic Structures on 3-manifolds, I:\\
Deformation of acylindrical manifolds}
\author{William P. Thurston}
\date{{\today}  reconstruction of 1986 Annals paper}
\address{
Mathematics Department \\
University of California at Davis \\
Davis, CA 95616
}
\email{wpt@math.ucdavis.edu}
\begin{abstract}This is an eprint approximation to \cite{Th1}, 
which is the definitive form of this paper. This eprint is provided
for convenience only; the theorem numbering of this version is different,
and not all corrections are present, so any reference or quotation
should refer to the published form. Parts II
and III (\cite{Hype2} and \cite{Hype3}) of this series,
although accepted for publication, for many years have only existed in
preprint form; they will also be made available as eprints.
\end{abstract}
\maketitle
\setcounter{section}{-1}
\section {Introduction}

This is the first in a series of papers dealing with the conjecture that all
compact $3$-manifolds admit canonical decompositions into geometric pieces.
This conjecture will be discussed in detail in part IV.
Here is an easily stated special case, in which no decomposition is
necessary:
\begin{conjecture}[Indecomposable Implies Geometric]
\label{indecomposable implies geometric}
Let $M^{3}$ be a closed, prime, atoroidal 3-manifold.
Then $M^{3}$ admits a locally homogeneous Riemannian metric,
{\it i.e.,}
a Riemannian metric such that for any two points $p$ and $q$ there is
an isometry from a neighborhood of $p$ to a neighborhood of $q$
carrying $p$ to $q$.
\end{conjecture}
\par
The main result of this series of papers will be to prove conjecture
\ref{indecomposable implies geometric}
and generalizations for a large class of manifolds.
\begin{theorem}[Atoroidal Haken Is Hyperbolic]
\label{atoroidal Haken is hyperbolic}
Let $M^{3}$ be a closed, atoroidal Haken manifold.
Then $M^{3}$ admits a hyperbolic structure,
i.e.,
a Riemannian metric with all sectional curvatures equal to -1.
\end{theorem}
\par
A hyperbolic structure is, of course, locally homogeneous: it is locally
isometric to hyperbolic 3-space $\hy^{3}$.
\par
Conjecture \ref{indecomposable implies geometric}
is far from being proven in general, but it seems well-supported by a great
many examples.
A special case is
\begin{conjecture}[Poincare]\label{Poincare}
Every simply-connected closed three manifold is homeomorphic
to the $3$-sphere.
\end{conjecture}
This would easily follow from
\ref{indecomposable implies geometric}
because it is easily seen that the only 3-dimensional simply-connected
locally homogeneous space is the three sphere (up to diffeomorphism).
The Poincar\'e conjecture has been the goal of many efforts by good
mathematicians, which
on the whole have not been very fruitful in shedding light on 3-manifolds.
It has probably been too narrow a goal;
I hope that the geometrization conjecture, which applies to all 3-manifolds
and hence has an abundance of examples to be tested against, will prove more
productive in the long run toward understanding 3-manifolds.
\par
We will also prove generalizations of
\ref{atoroidal Haken is hyperbolic}
to manifolds which are not closed, as well as to orbifolds or V-manifolds
(which are related to discrete groups of isometries which do not act freely).
\par
There has been a rather long gap in time since the time I first
formulated the
geometrization conjecture and announced and discussed its proof for
Haken manifolds in various forums and the present.
In the meantime, the theorem
has become widely known and accepted,
although relatively few people have learned the proof.
Part of the difficulty with writing about it has been that there is a lot
of background material which feeds into the proof at least in the
informal sense, if not in the logical sense.
This background material is nice in itself --- it seems to be an intrinsic
part of the geometric theory of three-manifolds --- but it is daunting because
of its quantity.
Much but not all of this background material is informally published in
\refto {Th1}, particularly chapters 8 and 9,
but little of it has been formally published.
This series of papers is organized so that the most important and
least published ideas come first.
Later parts will present material tending more to
overlap with
\refto {Th1},
but the presentation will be more careful.
\par
Almost all proofs having to do with Haken manifolds take place by
induction --- by cutting the
manifold along incompressible surfaces into simpler
manifolds (with boundary), analyzing the resulting pieces,
and reassembling.
The present proof is no exception.
The hard part of the proof is the geometric analysis of all possible
hyperbolic structures on a manifold with boundary in order to find limits to the
way the geometry is distorted as the hyperbolic structure is varied.
\par
The necessary information will be obtained in parts I, II and III.
Part II will also construct hyperbolic structures on atoroidal 3-manifolds
which fiber over a circle.
This material has been exposed in
\refto {Sul2};
Sullivan introduced several important simplifications to my original
unpublished approach, and in part II Sullivan's exposition
will be further simplified and generalized.
Hyperbolic structures on atoroidal Haken manifolds will be assembled
in part IV.
The other seven kinds of 3-dimensional geometric structures will also
be explained here, to give a conjectural overall picture of 3-manifolds.
\par
Part V will develop some results on geodesic laminations on surfaces and in
3-manifolds, and some material on homeomorphisms of surfaces.
Some of this material is closely related to my preprint
\refto {Th2}
which was expanded and exposed in
\refto {FLP}.
Other parts of the material
are developed in
\refto {Th1}
\par
Part VI will develop some necessary
theory of geometrically tame hyperbolic manifolds.
The needed theory was first developed in
\refto {Th1};
about a year ago, Francis Bonahon
\cite{Bon}
has proven a very nice stronger result.
\par
Part VII will deal with orbifolds.
It will give a proof of the existence of a certain ``reflection groups''
crucially needed in a kind of doubling trick in part IV.
It will also a stronger result, that the geometrization conjecture
holds for orbifolds with nonempty singular locus.
The stronger result for orbifolds logically depends on the general result
for Haken 3-manifolds, so the weaker form must be proven first.
\par
During the long delay in publication of this material, several worthwhile
references have appeared.
Peter Scott has written an excellent expository article concerning
3-dimensional geometric structures
\cite{Scott:geometries},
and John Morgan
\cite{Mor}
has written a detailed outine of the proof of the geometrization
theorem for hyperbolic 3-manifolds.
I also
wrote an expository article concerning the geometrization conjecture
\cite{Th3}.
John Morgan and Peter Shalen
\cite{M-S}, \cite{M-S2}, \cite{M-S3}, \cite{M-S-intro}
have found an interesting algebraic proof
of the main result of the current paper, and the generalization in
part III.

\begin{figure}[htb]
\centering
\includegraphics{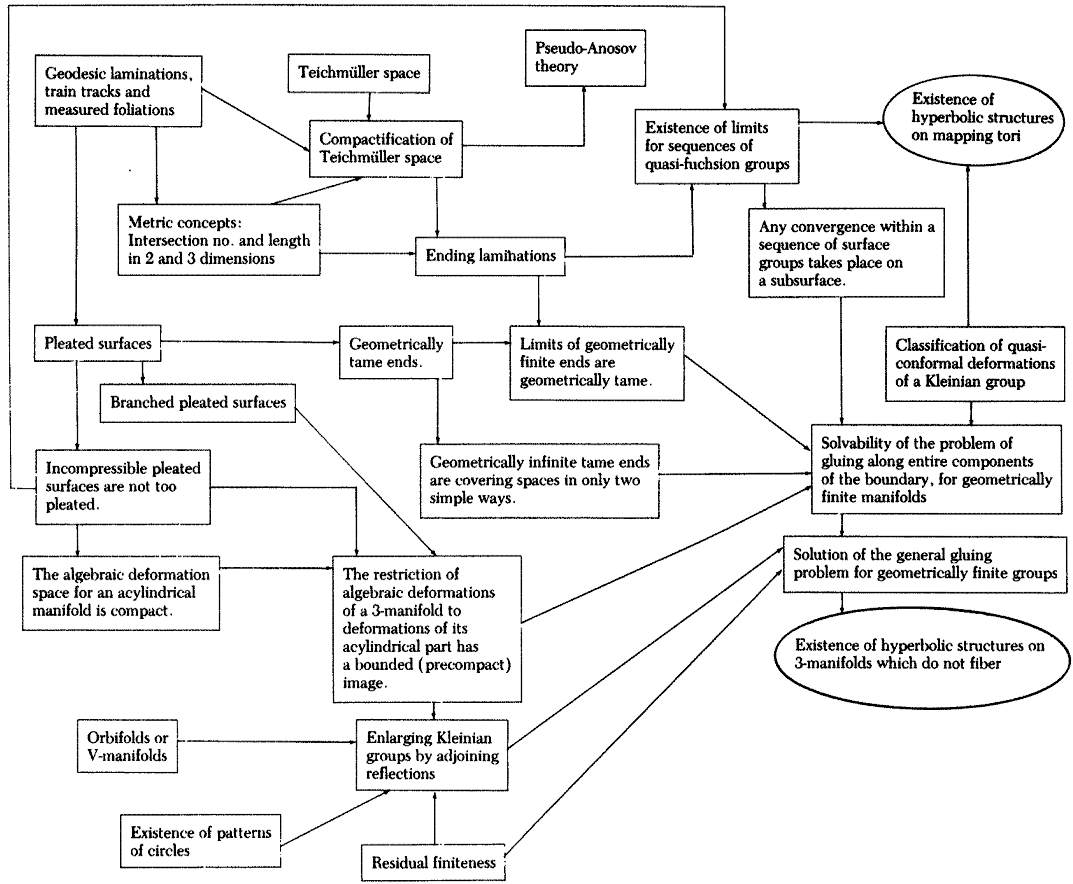}
\caption{This is a flowchart of some of the ideas which enter in the proof of
theorem \ref{atoroidal Haken is hyperbolic}.}
\label{Flowchart of ideas}
\end{figure}

Of somewhat indirect relevance is also
\refto{Th4}
in which I prove a result concerning conformal dynamical systems which closely
parallels many of the methods used in the current series of papers, although in
a much simpler form. A proof of this theorem has also been presented by Douady
and Hubbard,
\cite{D-H}.

I would like to thank George Francis for drawing most of the
illustrations.
\par
For purposes of reference, figure \ref{Flowchart of ideas}
describes the logical flow of ideas and results
which enter into the proof of the main theorem.
The boxes do not represent ideas of equal importance or difficulty,
and they could easily be regrouped or subdivided --- but as the
reader becomes familiar with the ideas of the proof, the chart
might help her or him to keep in mind the overall logical structure.

\section {Topologies}

For any $n$-manifold $M$, let $H ( M )$ denote the set of complete
hyperbolic $n$-manifolds $N$ equipped with a homotopy equivalence
$f : M \arrow N$.
Two elements $N$ and ${N'}$ are equivalent if there is an isometry
between them in the correct homotopy class,
that is, in the class that makes a homotopy commutative diagram.
$H ( M )$ is the same as the set of conjugacy classes of discrete
faithful representations of $\pi_{1} ( M )$ in the group of isometries
of hyperbolic $n$-space $\hy^{n}$.
\par
There are at least three distinct useful topologies for $H ( M )$;
glossing over the distinctions is a common source of error in arguments
involving varying hyperbolic structures.
We will define and briefly discuss in turn the algebraic topology, the geometric
topology, and the quasi-isometric topology.
\par
When $H ( M )$ is thought of in terms of representations of the fundamental
group of $M$,
an obvious topology comes to mind, namely the
topology of convergence of representations on finite generating sets.
This can be alternately phrased
as the compact-open topology on homomorphisms between the two groups.
The set $H ( M )$ inherits a topology as a quotient
space under the action of $\Isom ( \hy^{n} )$ by conjugacy on the
space of representations.
This topology is the
{\it algebraic topology}.
$H ( M )$ with this topology will be denoted $AH ( M )$.
In the case $M$ is a closed 2-manifold, $AH ( M )$ is commonly
called the
{\it Teichm\"uller space of $M$}.
\begin{proposition}[AH(Non-Elementary) Is Hausdorff]\label{AH(non-elementary) is Hausdorff}
If $\pi_{1} ( M )$ does not contain an abelian subgroup with
finite index, then $AH ( M )$ is Hausdorff.
\end{proposition}
\par
A discrete
group of isometries of hyperbolic space which admits an abelian subgroup
with finite index is called an
{\it elementary}
group.
\begin{proof}
Suppose that $\pi_{1} ( M )$ does not have an abelian subgroup of
finite index and that $N_{1}$ and $N_{2}$ are elements of $AH ( M)$
such that every neighborhood of $N_{1}$ intersects every neighborhood
of $N_{2}$.
We will show that $N_{1}  =  N_{2}$.
\par
We translate to the context of representations of groups.
Suppose that $M$ does not contain an abelian subgroup of finite index.
Represent $N_{i}$ as $\hy^{n} / \rho_{i} ( \pi_{M} )$.
Let $g_{1} \dots g_{k}$ be generators for the fundamental
group of $M$.
For any representation $\tau$,
the translation distance for the $i$th generator
is a function $T_{\tau} ( g_{i} )(x)  =  d ( x ,  \tau ( g_{i} ) ( x ))$
on hyperbolic space.
The translation distance of $g_i$
goes to infinity with distance from the axis of $\tau ( g_{i} )$ if
$\tau_{g_i}$ is hyperbolic. Otherwise, $\tau_{g_i}$
has a unique fixed pointed at infinity,
and $T_\tau(g_i)$ tends to infinity away from this fixed point.
\par
It cannot happen
that all generators for a discrete faithful representation $\tau$
have a common fixed point at infinity, for otherwise the entire group would
be a discrete subgroup of the group of similarities of Euclidean space, which
automatically has an abelian subgroup of finite index.
(This is a standard part of the theory of the thick-thin decomposition,
or Margulis decomposition, of a hyperbolic manifold).
Therefore, the total translation distance
$$d_{\tau}  =  \sum_{i} T_{\tau} ( g_{i} )$$
is a proper function, that is, it goes to infinity near the sphere at
infinity.
Also, $d_{\tau}$ is strictly convex, so that $d_{\tau} ^{-1} ( (0,  K] )$
is a convex ball in $\hy^{n}$ for any sufficiently large $K$.
\par
Suppose that a sequence of representations $\phi_{i}$ converges to
$\rho_{1}$, while a conjugate sequence $\psi_{i}$ converges to
$\rho_{2}$.
Then the sets $d_{\phi_{i}} ^{-1} ( ( 0 ,  K ])$ must converge to
$d_{\rho_{1}} ^{-1} ( ( 0 ,  K ])$ in the Hausdorff
topology,
and similarly for $\psi_{i}$.
The conjugating isometries must carry these sets for the $\phi_{i}$
to the corresponding sets for the $\psi_{i}$.
This implies that the conjugating isometries remain in a compact subset
of the group of isometries of hyperbolic space, so they have a convergent
subsequence,
whose limit conjugates $\rho_{1}$ to $\rho_{2}$.
\end{proof} 
\par
The case that $M$ is a closed surface has been analyzed classically.
Fricke proved that for a closed surface of negative Euler characteristic,
$AH ( M )$, which may be identified with the Teichm\"uller space
$T ( M )$ for $M$, is homeomorphic to Euclidean space of
dimension $3 |\chi (M)| $.
When $M$ is a surface with non-empty boundary, $AH ( M )$ is not connected.
The explanation is that a homotopy equivalence between
two surfaces with boundary is not generally homotopic
to a homeomorphism.
The components of $AH ( M)$ are in one-to-one correspondence
with homotopy equivalences of
$M$ to surfaces, up to homeomorphisms.
Each component of $AH ( M )$ is diffeomorphic to the product of
Euclidean space with a number of half-lines.
\par
A three-manifold $( M , \boundary M )$ is
{\it acylindrical}
if $\boundary M$ is incompressible and if every map
of the 2-dimensional cylinder $C  =  S^{1} \times I$,
$$f: ( C , \boundary C )  \arrow  ( M ,  \boundary M )$$
which takes the components of $\boundary C$ to non-trivial
homotopy classes in $\boundary M$ is homotopic into $\boundary M$.
The main result of this paper is
\begin{theorem}[AH(Acylindrical) Is Compact]\label{AH(acylindrical) is compact}
If $M$ is any compact acylindrical 3-manifold with boundary,
then $AH ( M )$ is compact.
\end{theorem}
The case when $\boundary M$ is empty follows from Mostow's theorem
\cite{Mos},
which asserts that $AH ( M )$ consists of at most one point.
The theorem is
false without the hypothesis that $M$ is acylindrical;
a simple example is $( surface  \times  I )$.

Another appealing and often useful topology for $H ( M )$
is the
{\it geometric topology},
denoted $GH ( M )$.
To define a neighborhood basis for $GH ( M )$, consider
$( N ,  f )$ be any element of $H ( M )$, where $f$
is a homotopy equivalence from $M$ to the hyperbolic manifold $N$.
A (small) neighborhood of $( N ,  f )$ in the geometric topology
is defined by the choice of a (large)
compact set $K  \subset  N$ and a (small) positive real number
$\epsilon$.
The neighborhood $\Nbhd_{K , \epsilon} ( N )$
consists of those elements $( {N'} ,  {f'} )$ of $H ( M )$
such that there exists a diffeomorphism $\phi$
of $K$ to a subset ${K'}  \subset  {N'}$
whose derivative is within $\epsilon$ of being
an isometry
at every point of $K$ and where $\phi$ has the correct homotopy class.
Explicitly,
$ (f')\inverse \circ i' \circ \phi$ must be homotopic to $f\inverse \circ i$,
where $i: K \subset N$ and $i': K' \subset N'$ are the inclusions
and $f\inverse$ and $(f')\inverse$ are homotopy inverses of $f$ and $f'$.
To measure the deviation of a linear map $L$
from being an isometry, one can use
the quantity
$$\sup_{V \ne 0} {\left | \log {\left ( \frac{\| L(V) \|}{\| V \|} \right )} \right |}$$

This geometric topology is not to be confused with
the geometric topology which can be
defined on other sets.
The most useful example
is the set of based hyperbolic manifolds, with no specification
of homotopy type or homotopy classes of maps.
The definition of the geometric topology on that set differs in that
the map $\phi$ must preserve base point, and its homotopy class
is not restricted.

A third topology on the set $H ( M )$ is the
{\it quasi-\-isometric topology}, denoted $QH ( M )$.
A small quasi-\-isometric neighborhood of an element
$f : M \arrow N$ of $H ( M )$ consists of those
${f'} : M \arrow {N'}  \in H ( M )$ for which there exists
a global diffeomorphism $g : N \arrow {N'}$ with
$g \circ f \homotopic {f'}$, whose derivative is uniformly
close to being an isometry.
\par
In general,
a diffeomorphism $g$ between two hyperbolic manifolds $P$ and $Q$
which has the property
that its derivative is within a uniform distance of an isometry is called
a
{\it quasi-\-isometry}.
More particularly, when the derivative of $g$ satisfies
$${\left | \log {\left ( \frac{\| dg ( V ) \|}{\| V \|} \right )} \right |}
 \le  K \quad\quad\text{for all tangent vectors $V$ to $P$,}$$
$g$ is called a
{\it K-quasi-\-isometry}.
When such a quasi-\-isometry exists, then $P$ and $Q$ are
{\it quasi-\-isometric}.
\par
These three topologies, the algebraic, geometric, and quasi-\-isometric,
are successively finer, so that the maps
$$QH ( M )  \arrow  GH ( M )  \arrow  AH ( M )$$
induced from the identity on $H ( M )$
are continuous.
This is an easy exercise from the definitions.
Since $AH ( M )$ is Hausdorff, it follows that the other two are also
Hausdorff.
\par
The inverse maps, in general, are discontinuous.
For example, it is easy
to see that in the geometric topology any non-compact hyperbolic
surface of finite area is a limit of hyperbolic surfaces with infinite area.
In the quasi-\-isometric topology, area is clearly a continuous function
(where the values of the area function are endowed with the topology making
$\infty$ an isolated point).
\par
The demonstration that $AH(M)  \arrow  GH ( M )$ is not continuous is
much trickier, and
requires at least three dimensions.
There exist examples for which $M$ is the product of a surface of genus 2
with an interval.
See
\refto {Th-K}
for one discussion of such an example.
\par
The sphere at infinity
divides into
two essentially different parts under the action of any discrete
group $\Gamma$ of hyperbolic isometries: the limit set $L_{\Gamma}$,
which is closed,
and its complement, the region of discontinuity $D_{\Gamma}$, which is open.
$D_{\Gamma}$ is the unique maximal open set where
$\Gamma$ acts properly discontinous, and
$L_{\Gamma}$ is the unique non-empty minimal closed subset invariant
for $\Gamma$ (except in certain cases when $\Gamma$ is elementary,
such as $\Gamma  =  \integers$).
\par
Since $\Gamma$ acts conformally on $D_{\Gamma}$, the quotient
space $D_{\Gamma} / \Gamma$ inherits a conformal structure.
Let us specialize to the case of dimension $n = 3$.
The {\it Teichm\"uller space of $( D_{\Gamma} ,  \Gamma )$} can be defined
as the space of all conformal structures on $D_{\Gamma}$ invariant by
$\Gamma$
which are quasiconformal to the initial conformal
structure, up to conformal equivalences between structures
which are equivariant with respect to $\Gamma$.
Except in the case of elementary groups, each such conformal structure
is represented by a unique $\Gamma$-invariant
hyperbolic metric on $D_{\Gamma}$, its Poincar\'e metric.
\par
The Teichm\"uller space $T ( D_{\Gamma} , \Gamma )$ is sometimes the same
as the Teichm\"uller space of the quotient surface $D_{\Gamma} / \Gamma$,
but not always.
Consider, for simplicity, the case
that $\Gamma$ has no elliptic or parabolic elements.
Then the Ahlfors finite area theorem says that
the quotient surface consists of a finite number of components,
each compact.
Any conformal structure on the quotient surface gives rise to an equivariant
conformal structure on $D_{\Gamma}$, and all equivariant conformal
structures are quasiconformally related.
This says that an element of the Teichm\"uller space of the quotient
surface defines an element of the Teichm\"uller space of
$(D_{\Gamma} ,  \Gamma )$.
When the components of $D_{\Gamma}$ are simply-connected, then any
equivariant homeomorphism of $D_{\Gamma}$ is equivariantly isotopic
to the identity, and the two Teichm\"uller spaces
are the same.
Otherwise, $T ( D_{\Gamma} ,  \Gamma )$ is the quotient of
$T ( D_{\Gamma} / \Gamma )$ by the group of isotopy classes of
quasiconformal homeomorphisms of $D_\Gamma$ which commute with the action
of $\Gamma$.

Any quasi-\-isometry between two complete hyperbolic n-manifolds lifts
to a quasi-\-isometry from $\hy^{n}$ to $\hy^{n}$.
Any such quasi-\-isometry extends continuously to the sphere at infinity,
where it induces a quasiconformal map.
(Indeed,
this fact was an important point in Mostow's proof that hyperbolic
structures are rigid.)
Therefore, if $N$ is any hyperbolic 3-manifold and ${N'}$ is
any other hyperbolic 3-manifold quasi-\-isometric with $N$, the conformal
structure on $D_{\pi_{1} ( {N'} )}$ defines an element
$${\rm conf} ( {N'} )  \in  T ( D_{\pi_{1} ( N )} , \pi_{1} ( N ) )$$
We shall be concerned with the case that each component of
$D_{\pi_{1} ( N )}$ is simply-connected, so the conformal
invariant is more directly thought of as element of the Teichm\"uller
space of the quotient surface.
\par
Here is a fundamental result, which has evolved through the work
of a number of people, beginning with Ahlfors and Bers, including
Mostow, Maskit and others, and put in a clean form through the
work of Sullivan.
\begin{theorem}[Quasiconformal Deformation Theorem]\label{quasiconformal deformation theorem}
Let $M$ be any complete hyperbolic 3-manifold with finitely generated
fundamental group.
The component $QH_{0} ( M )$ of $QH ( M )$ which contains $M$
consists of all $f : M \arrow {M'}$
such that $f$ is homotopic to a
quasi-\-isometry $\phi : M \arrow {M'}$.
The map
$${\rm conf} :  QH_{0} ( M )   \arrow   T ( D_{\pi_{1} ( M )} , \pi_{1} ( M ) )$$
is a homeomorphism.
\end{theorem}
\par
In particular, the quasi-\-isometric deformation space $QH_{0}$
always consist of a single point, or it is noncompact.
An important special case occurs when $M$ is compact or has
finite volume.
Then $D_{\Gamma}$ is empty, so $QH_{0} ( M )$ is a single
point; this case is the Mostow rigidity theorem.
Since $QH_{0}$ is often not compact, theorem
\ref{AH(acylindrical) is compact}
may be thought of as an existence theorem for many sequences
of elements of $H ( M )$ which converge in $AH ( M )$ but
not in $QH ( M )$.
\par
A complete description of the three spaces $AH ( M ),  GH(M),$
and $QH(M)$ is certainly not rigorously known, but here is a
conjectural image, of which certain features can be rigorously
proven.
Let us stick to the case that $M$ is a compact, acylindrical
manifold.
Then $H(M)$ is a hard-boiled egg.
The egg complete with shell is $AH(M)$; it appears to be homeomorphic
to a closed unit ball.
$GH(M)$ is obtained by thoroughly cracking the egg shell on a convenient
hard surface, such as one's skull.
Apparently no material is physically separated from the egg, but many
cracks are developed --- cracks are dense in the boundary --- and
at the same time, the material of the egg just inside the shell is
weakened, so that neighborhood systems of points on the boundary
become thinner.
Finally, $QH(M)$ has uncountably many components, which are obtained
by peeling off the shell and scattering the pieces all over.
Each component is homeomorphic to some Teichm\"uller space --- it
is parametrized by Euclidean space of some even dimension.
``Most'' of the components have dimension zero, for they describe
groups whose limit set is all of $S_{\infty}^{2}$.

\section {Ideal triangulations}

The plan for the proof of theorem \ref{AH(acylindrical) is compact} is to study
maps of a triangulation of a compact acylindrical manifold $M$ to hyperbolic
manifolds homotopy equivalent to $M$. In a sequence of such maps, we can study
how the shapes of simplices might degenerate. Hyperbolic simplices can
degenerate in only very special ways; qualitative information about the
geometry of the three-manifolds can be deduced from these shapes and how they
are assembled.

It is desirable to impose conditions on the maps we study so that properties of
the map indicate properties of the target manifold, and are not so much
dependent on the arbitrary choices we make in choosing a particular map in a
homotopy class. From the beginning, we abandon any thoughts of choosing the map
to be an embedding, and concentrate on more attainable conditions.

If $\tau$ is a triangulation of $M$, and if $f: M \arrow N$ is any map in the
homotopy class, then $f$ is homotopic rel the vertices of $\tau$ to a {\it
piecewise straight} map $s$, that is, having the property that when $s$ is
lifted to a map $\tilde {s}$ between universal covers, the image of each
simplex is the convex hull of the image of its vertices. It is not important
exactly how each simplex maps to its image, but one way to choose an explicit
map is by way of the Lorenz model for hyperbolic space. The lift $\tilde {f}$
of $f$ to a map between universal covers sends the vertices of any simplex to
points which in the Lorenz model $\euclidean^{3,1}$ lie on $\hy^{3}$ regarded
as a sheet of a hyperboloid. The map on vertices extends to a unique map $S$ of
$\tilde {M}$ to $\euclidean^{3,1}$ which is linear on each simplex; when this
is projected back to the hyperboloid, it gives a map $\tilde {s}$ to $\hy^{3} 
=  \tilde {N}$. Since $\tilde {s}$ is canonically defined, it is equivariant,
so it projects to a map to $N$. $S$ is equivariantly homotopic to $\tilde {f}$
rel vertices, by a linear homotopy. The homotopy can be projected back to
$\hy^{3}$ and thence to $N$, giving a homotopy of $f$ to $s$ rel vertices.

It is clear that any other piecewise-straight map is homotopic to
$s$ through maps such that each simplex has the same image.
This property also holds for the lifts of the maps to the universal covers.

The use of piecewise straight maps removes much of the arbitrary choice
involved in choosing a map in the homotopy class --- the choice which remains
is the choice of positions of the vertices, except for the irrelevant choice of
parametrization of images. To be precise, the choice of positions of vertices
is really made not in $N$, but in its universal cover $\tilde {N}  =  \hy^{3}$.
Any equivariant map of the vertices of the triangulation $\tilde {\tau}$ of
$\tilde {M}$ determines a piecewise straight map, up to trivial homotopy. If
the images of the vertices are merely chosen in $N$, the homotopy classes of
edges must still be determined.

To minimize choices, we will use a limiting case of piecewise-straight maps
where we choose positions of the vertices not in $\hy^{3}  =  \tilde {N}$, but
on the sphere at infinity $S_{\infty}^{2}$. By doing this, we immediately will
lose the map from $M$ to $N$. We will have at most a map from $M$ minus the
vertices of the triangulation to $N$. What we will gain is much simpler
variations in the shapes of images of simplices.

An {\it ideal $k$-simplex} in $\hy^{n}$ is the convex hull of a set of $k+1$
distinct points on $S_{\infty}^{n-1}$, subject to the condition that this
convex hull is $k$-dimensional. If the convex hull is lower dimensional, then
it is a {\it flattened ideal simplex}. An ideal simplex is homeomorphic to a
simplex minus its vertices. To choose an explicit parametrization, one can
choose a point in $\euclidean^{n,1}$ to represent each vertex. These
representatives lie on the light cone; the choice can be anywhere on the
appropriate ray. Once these choices are made, map the simplex to
$\euclidean^{n,1}$ by linear extension, and project the simplex minus its
vertices to $\hy^{n}$. This determines the parametrization up to precomposition
with projective diffeomorphisms of the simplex to itself.

Similarly, there is a map of a simplex minus its vertices to a flattened ideal
simplex which is well-determined up to precomposition with projective
diffeomorphisms.

\begin{figure}[htb]
\centering
\includegraphics{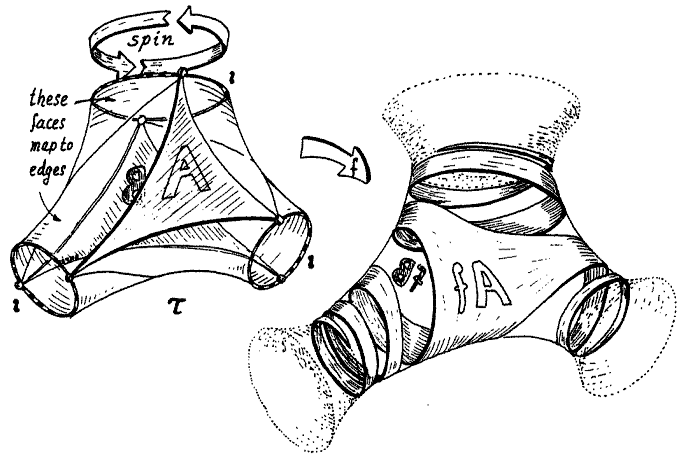}
\caption{An ideal simplicial
map for a pair of pants can be obtained by spinning a triangulation.
The initial triangulation has $8$ triangles, but only two of
them remain nondegenerate after spinning.}
\label{ideal simplicial pants}
\end{figure}

An {\it ideal simplicial map} of $\tau$ to $N$ is defined to be an equivariant
map from the universal cover $\tilde {\tau}$ to the completion $\hy^{3} \bigcup
S_{\infty}^{2}$ of $\tilde {N}$ such that each vertex of $\tau$ goes to
$S_{\infty}^{2}$, and the image of each simplex is the convex hull of the image
of its vertices. Any ideal simplicial map $r$ determines a special {\it
subcomplex at infinity} $\tilde {\iota} ( r )$ of $\tilde {\tau}$, consisting
of each simplex whose image is a single point at infinity. The image of
$\tilde\iota(r)$ in $\tau$ is denoted $\iota ( r )$. We shall also refer to
$\iota$ as the subcomplex at infinity. There is a map $D(r)$ defined
downstairs, but its domain is only $\tau - \iota$. Each simplex of $\tau$ minus
its intersection with $\iota$ maps either to an ideal simplex, a flattened
ideal simplex, or in the case its intersection with $\iota$ contains more than
vertices, it maps degenerately to a lower-dimensional simplex. The map $D(r)$
is a limit of ordinary piecewise-straight maps restricted to $\tau  -  \iota$.

\begin{figure}[htb]
\centering
\includegraphics{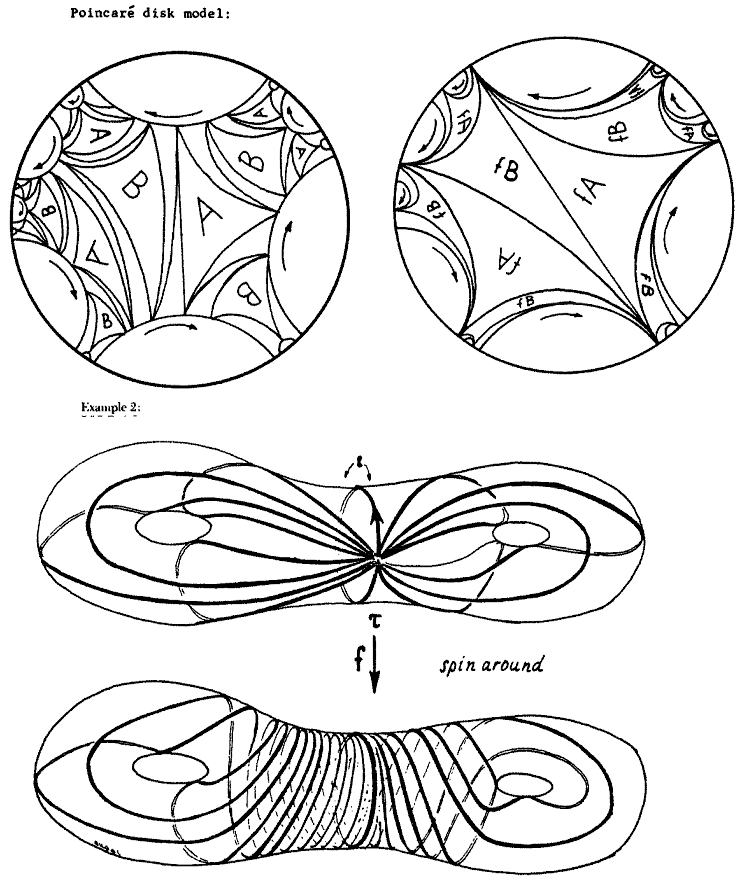}
\caption{Examples of ideal simplicial maps for surfaces. These can be
constructed starting with an ordinary piecewise straight map, then homotoping
the vertices through homotopy classes tending to infinity.}
\label{2-dimensional ideal simplicial maps}
\end{figure}

Figures \ref{ideal simplicial pants} and
\ref{2-dimensional ideal simplicial maps}
illustrate some examples in 2 dimensions.
These examples have been obtained from piecewise straight maps of
$\tau$ by spinning, that is, homotoping each vertex repeatedly around
a simple closed curve on the target surface and passing to the limit.
In example 1 on a pair of pants,
only 2 out of 8 two-simplices are mapped non-degenerately,
and in example 2 on a surface of genus 2, only 4 out of 6 are
non-degenerate.
These examples were chosen to be particularly nice, so that the
nondegenerate simplices are mapped disjointly and with positive
orientation.
In the general case for a surface of genus 2, each ideal triangle
would be
likely to be dense on the surface, although most of it would be
be extremely thin.

For the proof of theorem
\ref{AH(acylindrical) is compact},
we will study the limiting geometry of the map $r$ or $D(r)$,
and try to reconstruct a limiting manifold from this geometry.
Before trying to pass to limits, however, we should ask ourselves
whether we can even reconstruct $N$ from the local geometry of
$D(r)$, that is, from the shape data
--- the shapes of nondegenerate
images of simplices
minus $\iota$ --- and the gluing data ---
the information describing how the nondegenerate
images of simplices fit together near any common face in $\tilde {N}$.
In example 2 of figure
\ref{2-dimensional ideal simplicial maps},
this information does not suffice to determine the hyperbolic structure.
In this case, $M  -  \iota$ is not connected, and it is impossible to
reconstruct how the two halves are joined together.

We shall see that
the situation is much better in dimension 3.
A group $G$ is an {\it amalgam} if it can be described as
either (1) a free product $G = A \underset C \to * B$, or (2) an HNN amalgam
$G = A \underset C \to *$.  The group $C$ is the {\it amalgamating} group.
A representation as an amalgam is {\it trivial} if it is of type (1) and
$C$ includes into $A$ or $B$ as the entire group.

\begin{proposition}[Abelian Is Small In Acylindrical]\label{abelian is small in acylindrical}
The fundamental group of a compact aspherical atoroidal acylindrical
3-manifold with incompressible boundary
cannot be represented as a nontrivial amalgam
such that the amalgamating group is abelian or contains
an abelian subgroup with finite index.
\end{proposition}
\begin{proof} 

For any representation of a 3-manifold group as a nontrivial amalgam,
there is an incompressible surface whose fundamental group
is some subgroup of the amalgamating group, by a standard 3-manifold
argument. (Form a classifying space by
forming the product of a classifying space $X_C$ for $C$ with an interval,
and gluing its ends to the classifying space for $A$ or $A$ and $B$.
Map the three-manifold to the classifying space, and make the map
transverse to $X_C \times \{.5\}$.  The preimage of $X_C \times \{.5\}$
is a surface.  Apply the loop theorem repeatedly to construct homotopies
of the map which make this surface incompressible.)

If the amalgamating group is abelian or contains an abelian subgroup
of finite index, then the surface could only be a sphere, projective
plane, torus, Klein bottle, or disk, annulus or Moebius band.
The hypotheses rule out any such surface which cuts the 3-manifold
non-trivially.
\end{proof} 

\begin{proposition}[Iota Bounded By Abelian]\label{iota bounded by abelian}
If $S$ is any component of the boundary of the
regular neighborhood of the subcomplex
at infinity $\iota$ for any ideal simplicial homotopy equivalence
of a compact acylindrical 3-manifold $M$
to a hyperbolic 3-manifold $N$,
then the image of $\pi_{1} ( S )$ in $\pi_{1} ( N )$ contains
an abelian subgroup with finite index.
\end{proposition}

\begin{proof} 
Each component of $\tilde {\iota}$ maps to a single point on
$S_{\infty}^{2}$.
Therefore,
the fundamental group of any component of the boundary of $\iota$
acts on $\tilde {N}  =  \hy^{3}$ with a fixed point at infinity.
Any discrete group with this property contains an abelian subgroup
with finite index.
\end{proof} 

\begin{corollary}[Complement Of Iota Is Large]\label{complement of iota is large}
The fundamental group of some component of the complement of $\iota$
maps surjectively
to the fundamental group of $M$.
\end{corollary}

\begin{proof} 
For each component $S$
of the boundary of the regular neighborhood
of $\iota$, there is a representation of $\pi_{1} ( M )$ as either
an amalgamated free product or an HNN amalgam, where the amalgamating
group is the image of $\pi_{1} ( S )$ in $\pi_1(M)$.  By
\ref{iota bounded by abelian}, the amalgamating subgroup is
virtually abelian.  Thus, by \ref{abelian is small in acylindrical},
the amalgam is trivial.
This implies that $S$
must separate, and one of the two pieces must have the same image
for its fundamental group as $S$.
Throw away all the pieces of this type; what remains is a single component
the image of whose fundamental group is $\pi_{1} ( M )$.
\end{proof} 

Examples (rather artificial) can be easily constructed where the complement
of $\iota$ is not connected.
For example,
choose $\tau$ to be a triangulation of $M$ which admits an embedding of
a sphere in its 2-skeleton, and so that the ball bounded by the sphere
contains at least 1 vertex.
Choose the ideal simplicial map so that the vertices of each upstairs
sphere (in $\tilde {M}$) map to a single vertex on $S_{\infty}^{2}$,
while the vertices in the upstairs balls map to points distinct from
the image of its boundary and distinct from each other.
Somewhat more interesting
examples may be constructed using a solid torus whose fundamental group
injects in $\pi_{1} ( M )$.
In this case, the vertices on the bounding torus must be mapped to one
of the fixed points for the action of its fundamental on $S_{\infty}^{2}$.

We can now characterize exactly which subcomplexes $\iota$
of a triangulation
$\tau$ can be the subcomplex
at infinity for an ideal simplicial homotopy
equivalence.
Certainly $\iota$ must contain all the vertices of $\tau$, and it
must satisfy the condition
that if the boundary of any simplex of dimension 2 or more is in
$\iota$, then the simplex is in $\iota$.
(A simplex of dimension bigger than $1$ has a connected boundary, which
implies that if its boundary is in $\iota$, it maps upstairs to a single
point).
If $\iota$ also satisfies the necessary condition that the image of each the
fundamental group of each component
of $\iota$ in $\pi_{1} ( M )$ is virtually abelian,
then an ideal simplicial map is
easily constructed by mapping each component of $\tilde {\iota}$
to some point fixed by the action of its stabilizer in $\pi_{1} ( M )$
on $S_{\infty}^{2}$.

The final ideal simplicial map $r$ can easily have the property that each
image simplex by $D(r)$ is dense in $N$.
However, most of each image simplex consists of long, thin spikes;
because the spikes have small volume, we shall see that they have
little real consequence.

\begin{figure}[htb]
\centering
\includegraphics{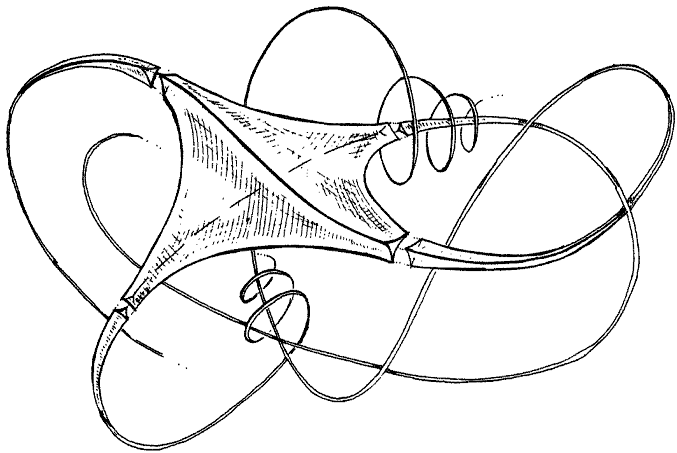}
\caption{The image of a simplex may be dense under an ideal simplicial map.}
\label{image may be dense}
\end{figure}


\section {Shapes and limits}
\label{Shapes and limits}

We will now analyze in detail how the shapes of ideal simplices
determine a hyperbolic manifold $N$, and what kind of limiting
behaviour they can have.
From this analysis, we will derive a quick proof of theorem
\ref{AH(acylindrical) is compact}
in the case that the manifold is compact; the theorem
is known as a consequence of the Mostow rigidity theorem in this case,
but the simple proof helps illustrate the idea.
\par
The remaining ingredient needed for the general case of
\ref{AH(acylindrical) is compact}
is some way to control the boundary of $M$.
We shall do this using an analysis of pleated surfaces,
to be introduced in the next section.
\par
For simplicity, we will stick to analyzing the case that $M$ is oriented.
The result for an unorientable manifold will be simple corollary of the
oriented case.
\par
All ideal 1-simplices in $\hy^{n}$ are congruent.
In fact, they are overly congruent: any ideal 1-simplex (line) has many
isometries taking it to any other ideal 1-simplex.
This is the reason that the local geometry of an ideal simplicial map does not
determine the target manifold in dimension 2.
In the case of an ideal 1-simplex in $\hy^{3}$, there is a two dimensional or
one complex dimensional group
of isometries stabilizing it.
\par
All ideal triangles in $\hy^{n}$ are congruent, since the group
of isometries of $\hy^{n}$ acts transitively on triples of
points.
In the case $n  =  3$, there is a unique orientation-preserving isometry
sending any ordered triple of points to any other.
\par
Ideal tetrahedra are not all congruent.
There is a one dimensional or one complex dimensional
family of different shapes.
The shape of an ordered simplex in oriented $\hy^{3}$
has a convenient parametrization
obtained by transforming by the unique
orientation-preserving isometry so the first three vertices of the simplex
at the points $0,  1$ and $\infty$ in the upper half-space model
Its congruence class is determined by the position of the last
vertex, which is at some point
$z  \in  \hat {\complexes}  -  \{ 0,  1,  \infty \}$
in the thrice-punctured sphere.
If $z$ is on the real line, then the ideal simplex is a flattened
simplex.
If the imaginary part of $z$ is positive, then the map of the
simplex preserves orientation, while if it is negative
it reverses orientation.
The invariant $z$ is the same as the cross-ratio of the four vertices
of the tetrahedron.
\par
There are 24 possible orderings of the vertices of a tetrahedron, 12
of them compatible with a given orientation.
How can the shape of an oriented but unordered simplex described?
\par
The description is simplified
by the fact that every ideal tetrahedron has a certain beautiful symmetry,
described by an action of $\zeemod 2 \times \zeemod 2$.
In fact, consider the common perpendicular to any opposite pair of
edges of an ideal simplex, that is, edges with no common ideal vertex.
Rotation of $\pi$ about such a perpendicular preserves the two edges
in question, which implies it preserves their set of endpoints.
Therefore, it acts as an isometry of the simplex to itself.
There are three pairs of opposite edges, so this construction gives
rise to three nontrivial isometries; together with the identity, they
form the group $\zeemod 2 \times \zeemod 2$.
The three axes intersect in a common point, where they are mutually
perpendicular.
\par
This group of symmetries of the ideal tetrahedron
is a normal subgroup of the alternating group for the four
vertices, and it preserves the invariant $z$.
The quotient group is $\zeemod 3$, which may be identified with the
group of even permutations of the
three pairs of opposite edges.
This means that the invariant $z$ is associated with a pair of opposite
edges.
\par
Here is a more direct way to describe the association of an
{\it edge invariant}
$z(e)$ to an edge $e$ of a congruence class of oriented ideal
simplices.
If $e$ is (temporarily) given an orientation, then in terms
of the orientation of the standard tetrahedron, the two triangles which
have $e$ as an edge can be distinguished as the clockwise face when
viewed from above and the counterclockwise face.
There is a unique orientation-preserving isometry $\phi$ of $\hy^{3}$ taking
the image of the clockwise face to the image of the counterclockwise face,
while preserving the two endpoints of $e$.
The edge invariant is the complex number whose modulus is the translation
distance of $\phi$, and whose argument is the angle of rotation of $\phi$.
Note that $\phi$ will rotate clockwise or counterclockwise about $e$,
depending on whether the simplex is mapped to preserve or reverse orientation.
If the opposite orientation of $e$ is chosen, then $\phi$ must be replaced
by $\phi \inverse$, but the complex number $z(e)$ remains the same.

\begin{figure}[htb]
\centering
\includegraphics{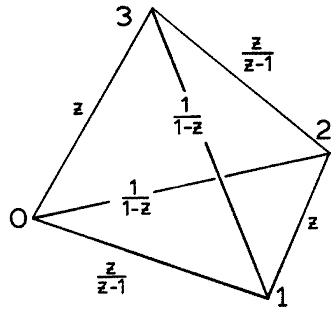}
\caption{The shape of an oriented ideal simplex is determined by any of its
edge invariants. The invariants for various edges are related as shown.}
\label{edge invariants}
\end{figure}

The edge invariant is consistent with the parametrization of an ordered
simplex:  if an ideal tetrahedron is arranged so the first three vertices
are at $0,  1$ and $\infty$ with the fourth vertex at $w$, then the
$w  =  z( \overline {0 \infty} )  =  z ( \overline {1 w })$.
The transformation taking the clockwise face for $\overline {0, \infty}$
to the counterclockwise face acts as multiplication by $w$ on the complex
plane.
\par
More generally, if one of the vertices is put at $\infty$ and its
then the invariant for the edge from any of the finite
points to $\infty$ is the ratio of the two adjacent sides of the
Euclidean triangle spanned by the three finite vertices;
here, a side of a triangle is interpreted as a vector, which
is identified with a complex number.
It is easy to deduce the formula
relating the various edge invariants to each other.
They are illustrated in figure~\ref{edge invariants}.

\begin{figure}[htb]
\centering
\includegraphics{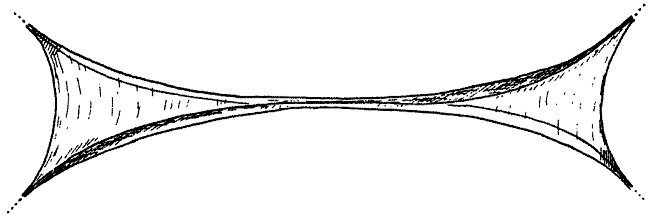}
\caption{As an ideal simplex changes, the only kind of degeneration which can happen
to its shape is that it stretches out into two parts which look like
ideal triangles, joined to each other by a long stringy part.
The edge invariants for the two edges which remain on one side
converge to $1$. The other edge invariants converge to $0$ or $infinity$,
depending on the orientation of the simplex.}
\label{stretched simplex}
\end{figure}

It is easy to understand the ways in which a single simplex can degenerate.
The illegal values for the edge invariant in $\Chat$ are $0$, $1$, and $\infty$.
Note that if one of the edge invariants tends toward any one of the
three illegal values, the other two edge invariants tend
toward the other two illegal values.
When this happens, the ideal simplex stretches apart into two
pieces that look nearly like ideal triangles, connected by a long thin
nearly 1-dimensional part (figure~\ref{stretched simplex}).

Now let $M$ be any compact three-manifold, and consider a sequence of elements
$\{ N_k \} \in \ah$, parametrized by ideal simplicial maps $f_k$ using domain
triangulation $\tau$. Let $\rho_k: \fund M \arrow \Isom \hy^3$ denote the
holonomy for $N_k$. We will pass to subsequences as necessary to try to get
convergence of as much of the data as possible; each subsequence will be
renamed like the entire sequence.  We first pass to a subsequence so that the
subcomplex $\iota$ at infinity is constant.  Next, we pass to a subsequence so
that the edge invariant of each image ideal simplex or flattened ideal simplex
converges in $\Chat$. If the limiting edge invariants are all in $\tps$, then
the sequence of representations converge up to conjugation, hence there is a
limit in $\ah$.

In general, some of the edge invariants may very well go to illegal values. We
will define a ``submanifold'' $G \subset M$ so that the sequence of
representations for the fundamental group of each component of $G$
automatically converges, possibly after taking a subsequence and conjugating.
For this, it is sufficient that $f_k$ restricted to $G$ is homotopic to a map
which stretches the metric of $G$ by only a bounded amount, since the set of
isometries which move a given point by a bounded amount is compact.

\begin{figure}[htb]
\centering
\includegraphics{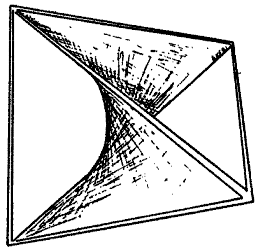}
\caption{A twisted square like this is adjoined to the 2-complex $B$ in each
$3$-simplex (of $\tau$) whose edge invariants tend toward the three illegal
values $0$, $1$, and $\infty$. The edges of the square are on the four
``stretched'' edges with invariants tending to $0$ and $\infty$, and the square
separates the two edges with invariants are tending to $1$.}
\label{twisted square}
\end{figure}

$G$ is obtained from $M$ by removing a neighborhood of a ``bad'' 2-complex
$B$.  $B$ is defined to consist of all $1$-cells (not because they are
bad, only because it simplifies the picture to discard them),
all degenerate $2$-cells, together with a twisted square embedded
in any $3$-simplex whose edge invariants are converging to illegal
values.  The twisted square is embedded in such a way that its four
edges are the four edges with edge invariants tending toward $0$
or $\infty$.  It separates the two edges whose edge invariants are
tending toward $1$ (figure \ref{twisted square}).

\begin{figure}[htb]
\centering
\includegraphics{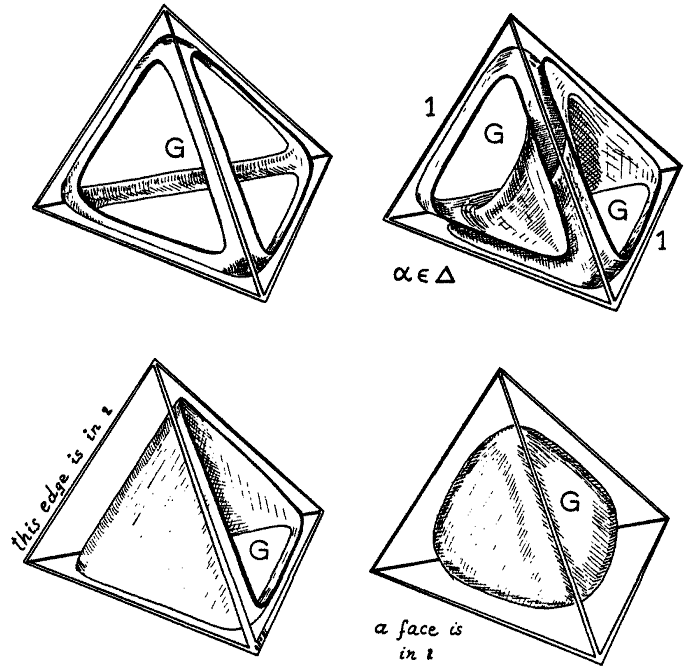}
\caption{These are the four ways which $G$ can intersect a $3$-simplex
$\alpha$, up to symmetries.  There are either $0$, $2$, or $4$ degenerate
simplices on the boundary of $\alpha$, depending on the intersection of $\iota$
with $\alpha$.  (Note that if $\alpha \cap \iota$ contains as much as edges or
a $2$-simplex, then all the $2$-simplices of $\boundary \alpha$ are
degenerate.)}
\label{G intersect simplices}
\end{figure}

Denote by $\Delta$ the set of all non-degenerate $3$-simplices with
edge invariants tending to illegal values.  The intersection of
$G$ with any 3-simplex $\alpha$ is determined by $\alpha \cap \iota$
and by whether or not $\alpha \in \Delta$.  There are four cases, illustrated
in figure \ref{G intersect simplices}.

It is clear from the definition that the sequence of representations restricted
to the fundamental group of a component of $G$ remains in a compact set.

Let $\kappa$ denote the frontier of $G$, that is,
$$ \kappa \cong \boundary G - \boundary M.$$
The strategy for proving theorem \ref{AH(acylindrical) is compact} will be to
show that $\fund G$ is large in $\fund M$; we will do that by showing that
$\kappa$ is in a certain sense ``small'' in $M$. The analysis is related to a
certain quotient map $p: \kappa \arrow \gamma$, where $\gamma$ is a graph or
$1$-complex.

\begin{figure}[htb]
\begin{centering}
\includegraphics{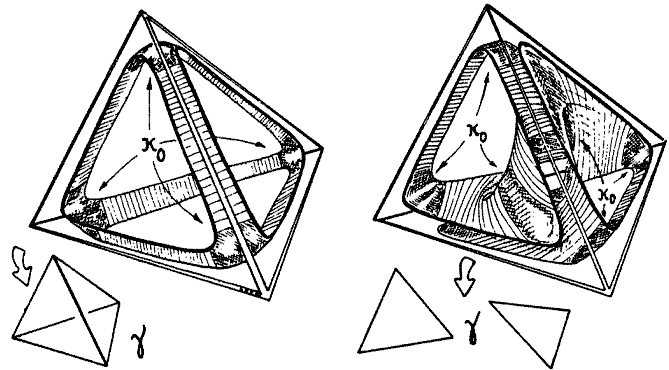}
\end{centering}
\caption{A quotient graph $\gamma$ of the internal boundary $\kappa$ of the
good submanifold $G$ is constructed by defining an equivalence relation in its
intersection with each $3$-simplex.}
\label{local decomposition of kappa}
\end{figure}

The quotient graph will be determined by defining equivalence relation
or decomposition of $\kappa$.

The first step is to construct a subsurface with
boundary $\kappa_0 \subset \kappa$ which consists of a the
intersection of a regular neighborhood of $\iota$ with $\kappa$,
together with a diagonal band on each side of the twisted squares of
$B$, as illustrated in figure \ref{local decomposition of kappa}.
Each diagonal band ``indicates'' vertices of ideal simplices which are
converging together, from the point of view of its component
of $G \cap$ the simplex. Each component
of $\kappa_0$ is to be collapsed to a point to form the vertices of
$\gamma$.

The second step is to form a foliation $F$ of
$\kappa - \kappa_0$.  In fact, any component of the
intersection of $\kappa-\kappa_0$ with a 3-simplex is a rectangle,
two of whose edges are on the boundary of the simplex and near two edges of
the simplex.  Let $F$ be the trivial foliation of each such rectangle
transverse to these two edges.  The leaf space of $F$ within the rectangle
is parametrized by the unit interval $[0, 1]$.
Arrange so that the foliations in neighboring rectangles
are joined so that this parametrization is either preserved or reversed.

The leaves of $F$ are compact, so each leaf is a circle if $M$ is closed,
or possibly an interval if $M$ has boundary.  The components of
$\kappa_0$ together with the leaves of $F$ define an equivalence relation,
whose quotient is the graph $\gamma$.

We will need the following, due to Chuckrow \cite{Chu}:
\begin{proposition}[Non-Elementary Limit Is Discrete]\label{non-elementary limit is discrete}
Let $\Gamma$ be any group, and suppose that $\sigma_k: \Gamma \arrow \Isom(\hy^n)$
is a sequence of discrete faithful representations which converges to a
a representation $\sigma_\infty$.  Then either $\sigma_\infty(\Gamma)$
is discrete and $\sigma_\infty$ is
faithful, or $\Gamma$ contains an abelian subgroup
with finite index.
\end{proposition}
\begin{proof}  For any $K > 0$ there is an $\epsilon > 0$ such that any
discrete group of isometries of $\hy^n$ which is generated by elements
$\gamma_1$, \dots, $\gamma_l$ having the property that for some point $x_0$,
$d(x_0, \gamma_1(x_0)) < \epsilon$ and for all $i>1$, $d(x_0, \gamma_i(x_0))<K$
has an abelian subgroup with finite index. ({\it Cf\/}.~\refto{Th1}, pp. ???
for one discussion).  To prove \ref{non-elementary limit is discrete}, let $K$
be greater than the distance that $\sigma_\infty$ of the generators of $\Gamma$
moves some $x_0$. If there is no element of $\Gamma$ such that $d(\sigma_\infty
\gamma x_0, x_0) < \epsilon$, then $\sigma_\infty \Gamma$ is discrete. 
Otherwise $\Gamma$ contains an abelian subgroup with finite index.
\end{proof}

\begin{proposition}[Gamma Locally Abelian]\label{gamma locally abelian}
For each $0$-cell or $1$-cell $\beta$ of the graph $\gamma$,
$\fund {p \inverse(\beta)}$ is abelian.
\end{proposition}
\begin{proof} If $\beta$ is an edge of $\gamma$, then $p \inverse (\beta)$
is a rectangle or a cylinder and in either case has an abelian fundamental
group.

Since $\kappa$ is in the good set, the sequence of representations $\rho_k$
converge on $\kappa$ to a limiting representation $\rho_\infty$. The limiting
representation restricted to any component of $\kappa_0$ fixes a point on
$\sinfty$. This point is a limit of images of vertices under the ideal
simplicial maps.  If $\fund{p\inverse(\beta)}$ is not discrete, it is abelian
by the preceding proposition. On the other hand, any discrete group of
orientation preserving isometries of $\hy^3$ which fixes a point at infinity is
abelian.
\end{proof}

It is likely that for some of the cells $\beta$ of $\gamma$, the group $\rho_k
\fund{p\inverse(\beta)}$ trivial. Let $\gamma_0$ be the union of cells whose
group is not trivial.

\begin{proposition}[Gamma0 Has Abelian Pieces]\label{gamma0 has abelian pieces} 
For each component $c$ of $\gamma_0$, the image
of $\fund{p\inverse(c)}$ in $\fund M$
is abelian.
\end{proposition}

\begin{proof}
This will follow from Proposition \ref{gamma locally abelian}
together with the following lemma:

\begin{lemma}[Commutation Transitive]\label{commutation transitive}

(i). If $x$, $y$ and $z$ are any three orientation preserving
isometries of infinite order in $\hy^3$
such that $x$ commutes with $y$ and $y$ commutes with $z$, then also
$x$ commutes with $z$.

(ii). If $x$ and $y$ are two orientation preserving isometries of infiite
order in $\hy^3$ such that $y\inverse x y$ commutes with $x$, and if
$x$ and $y$ generate a discrete group,
then $y$ commutes with $x$.
\end{lemma}

\begin{proof}

Orientation preserving isometries of infinite order commute iff they have the
same axis or the same parabolic fixed point at $\infty$.  Part (i) immediately
follows.

Under the hypothesis of part (ii), if $x$ is hyperbolic, then $y$ must also be
hyperbolic with the same axis, so they commute.  If $x$ is parabolic, then $y$
must fix the parabolic fixed point of $x$.  The hypothesis of discreteness
implies that $y$ cannot be hyperbolic; therefore it is parabolic, with the same
fixed point.

\end{proof}

The image of $\fund{p\inverse(c)}$ in $\fund M$ is built up by amalgamated free
products and HNN amalgams from abelian subgroups.   By the lemma, the entire
group is abelian.

\section{The deformation space for closed manifolds is compact}
\label{The deformation space for closed manifolds is compact}

We will now specialize and complete the proof of theorem
\ref{AH(acylindrical) is compact} in the case that $M$ is a
closed manifold or a manifold
whose boundary consists of tori.  This is equivalent to the condition that
$N_k$ have finite volume.  Even though this case follows from Mostow's rigiditry
theorem, it is worth presenting in order
to motivate the more technical version needed to handle the case when $M$
has boundary.

Here is a quick synopsis of the proof.  The geometry of the maps from
$M$ to $N_k$ (pinned down via ideal simplicial maps) remains bounded unless
something is stretched a great deal.  In terms of the simplices,
$M$ decomposes into pieces where the geometry of the map remains
bounded, joined by other parts of $M$ where the geometry is going to infinity.
The interface is a surface which maps with small area, and which turns out
to be made up of pieces whose image fundamental group in $M$ is
abelian.  By topological considerations, such surfaces cannot separate
an atoroidal manifold in an essential way.
Thus, one of the pieces where the geometry of the map remains bounded
carries all of the fundamental group of $M$.

Since $\boundary M$ consists of tori,
we can arrange that $\boundary M \subset \iota$,
by mapping all simplices of $\boundary \tilde M$ to the point at infinity
corresponding to the given cusp.  This implies that $\kappa$ will be
a closed surface.

We can make an isotopy of the surface $\kappa$ in $M$ so that its image
in $N_k$ by the ideal simplicial map has small area, tending to 0 as
$k \arrow \infty$ .  This can be accomplished
by isotoping the intersection of $\kappa$ with any $2$-simplex so that
it is near
near the boundary of the $2$-simplex, and
then extending individually to each $3$-simplex.
Recall that components of $\kappa - \kappa_0$ intersected with a
$3$-simplex are foliated rectangles.  We arrange that these
rectangles are mapped to long thin strips in $N$, whee the leaves of the
foliation give the short direction.  In fact, we make the lengths of the
rectangles go to infinity with $k$, while the widths go to $0$.

For any edge $e$ of $\gamma - \gamma_0$,
the cylinder $p\inverse(e) \subset \kappa$
is made up of a fixed number of thin strips, glued together edge to edge.
These edges are the intersections of the cylinder with $2$-simplices
of $\tau$.  Since the fundamental group of the cylinder has a trivial
representation in this case, $f_k$ restricted to the cylinder can be
lifted to $\hy^3$.  Since the cylinder lies
in the good set, the amount of shearing of one of these
strips in the cylinder with its neighbors
is bounded as $k \arrow \infty$.
Since the strips are becoming infinitely long, there is a short
cross sectional circle.
Such a circle bounds a disk of small area.

Construct a manifold $G_1$ by attaching $2$-handles to $G$ along
the core curve of each of these trivial cylinders.  The inclusion
$G\subset M$ extends to a homotopically unique map $g_1: G_1 \arrow M$,
since the core curves are trivial in $\fund M$.  The composed map
$f_k \circ g_1$ into $N_k$ is homotopic to a map such that the
image of $\kappa_1 = \boundary G_1$ has small area.  Note that
the image of the fundamental group of any component of $\kappa_1$
in $\fund M$ is abelian, in virtue of \ref{gamma0 has abelian pieces}.

\def\mass{\mathop{\rm mass}\nolimits}
If $C$ is any piecewise-differentiable singular $3$-chain
in a Riemannian manifold $P$, let $\degree_C: P \arrow \reals$
be the step function which gives the degree of $C$ at a point in $P$.
It is defined almost everywhere.  The {\it mass} of $C$ is defined as
$\mass ( C) = \int_P |\degree_C| \, dV$.  The mass of a $2$-chain
is defined similarly, using the degree of $C$ on each tangent $2$-plane
in $P$ and integrating with respect to $2$-dimensional Lebesgue measure.
For a generic $2$-chain, the mass is just the area of its image.

\begin{proposition}[Small Abelian 2-Cycle Bounds Small 3-Chain]\label{small abelian 2-cycle bounds small 3-chain}
Let $P$ be any hyperbolic $3$-manifold with abelian fundamental group.
Then every $2$-cycle $Z$ in $P$ bounds a $3$-chain $C$ in $P$,
possibly with non-compact support, satisfying $\mass(C) \le \mass(Z)$.
The function $\degree_C$ is uniquely determined by $Z$.
\end{proposition}
\begin{proof}
There are two cases --- the fundamental group of $P$ can be either parabolic
or hyperbolic.
Consider first the hyperbolic case.  Let $(r, x, \theta)$ be cylindrical
coordinates about a line in $\hy^3$, so $r$ is the distance from the line,
and $x$ is the height of the projection to the line.  The hyperbolic
metric is
$$ ds^2 = dr^2 + \cosh^2 r dx^2 + \sinh^2 r d\theta^2. $$
The $2$-form $\alpha = \cosh r \sinh r dx d\theta$, which is the area form
for equidistant surfaces about the line, has norm $1$.
Its exterior derivative
$$ d\alpha = (\cosh^2 r + \sinh^2 r) dr dx d\theta $$
is greater than the volume form
$$ dV = \cosh r \sinh r dr dx d\theta.$$
The form $\alpha$ descends to $P$.

If $Z$ is any $2$-cycle, $Z$ is the boundary of a compactly-supported
$3$-chain whose degree function is uniquely determined.
It suffices to prove the inequality for the positive and negative parts of
$\degree_C$, after perturbing $C$ to make $\boundary C$ transverse to itself.
Since we can change the sign of the negative part, we need consider
only the case $\degree_C \ge 0$.  We have
$$ \mass(C) = \int_C \, dV < \int_C \, d\alpha = \int_Z \, \alpha \le \mass(Z).$$
The case that $P$ is parabolic is similar,  except that $H_2(P)$ may not be zero; howerver,
each $2$-cycle bounds a unique chain of finite volume.  For this case, take as
$\alpha$ the area form for a family of horospheres.  Then $d\alpha=2dV$, and the proof goes through.
\end{proof}
\medskip

For each component $C$ of $\kappa_1$, the group $g_{1*} \fund C$ in $\fund M$ is
abelian --- either $1$, $\integers$, or $\integers^2$.  In the first two cases,
choose some $3$-manifold $M_C$ with $\boundary M_C = C$ such that the homomorphism extends to $\fund {M_C}$.  The map of $C$ to $M$ extends over $M_C$
in a homotopically unique way.  In the third case,
define $M_C = C \times [0, 1]$, and glue it on to $C$ at $C\times 0$.
The map $g_1$ on $C$ extends over $M_C$ taking $C \times 1$ to $\boundary M$.

Let $G_2$ be the manifold obtained from $G_1$ by attaching the $M_C$'s, and let
$g_2: (G_2, \boundary G_2) \arrow (M, \boundary M)$ be the extension of $g_1$
over the $M_C$'s. Clearly each component of $G_2$ contains a unique component
of $G$ and the corresponding components of $G_2$ and $G$ have fundamental
groups with the same image in $\fund M$. Thus, if $G_2^0$ is a component of
$G_2$ then, after conjugations, the representations $\rho_k \circ g_2: \fund
G_2^0 \arrow \Isom (\hy^3)$ converge.

\begin{proposition}[Degree G Is 1]\label{degree g is 1}
The degree of $g_2: (G_2, \boundary G_2) \arrow (M, \boundary M)$ is $1$.
\end{proposition}
\begin{proof} The degree can be computed by composing with the map
$f_k: (M, \boundary M) \arrow (N, {\rm cusps})$.  The map $f_k$, as an
ideal simplicial map, is not continuous on all of $M$, but by modifying
$f$ in a small neighborhood of $\iota$, it can be approximated by a
continuous map which is a homotopy equivalence.  The modification
takes place on a portion of $M$ whose image has a small mass.
The degree of $f_k$ can be defined at almost every
point in $N_k$ by adding up the degrees of the $3$-simplices which hit it.
It follows that the local degree of $f_k$ is $1$ almost everywhere.

The singular chain defined by $g_2$ is obtained from that defined by $f_k$
by a sequence of moves, each of which affects only a small volume:
The mass of the image of $B$ (the complement of $G$) is small; $G_1$ was
obtained by adding $2$-handles with small mass; and finally, $G_2$
was obtained from $G_1$ by adding 3-chains of small mass.
Thus, the total mass of the chain defined by $g_1$ is close to the volume
of $N_k$.  Since this mass must be an integral multiple of the
volume of $N_k$, and since
there is an {\it a priori} bound to the volume of a hyperbolic
manifold,  it follows that its degree must be 1.
\end{proof}

This proposition tells us that $G_2$ must ``capture'' most of the
fundamental group of $M$, so the proof is almost finished.  In fact,
if $(P, \boundary P) \arrow (Q, \boundary Q)$
is a map of degree $k$ between connected $n$-manifolds, the image of $\fund P$
in $\fund Q$ has index at most $k$: otherwise, the map would lift to
a covering space of degree (possibly infinite) greater than $k$,
so the composed map from $P$ to $Q$ would have degree greater than
$k$.  The proof in this special case ($\boundary M$ consists of tori)
of theorem \ref{AH(acylindrical) is compact} is finished, except for
the possibility that $G$ and hence $G_2$ is not connected.

Even if $G_2$ is not connected, at least one component must map to
$M$ with positive degree.  The easy case of theorem
\ref{AH(acylindrical) is compact} will be completed by the following:

\begin{lemma}[Finite Index Convergence Implies Convergence]\label{finite index convergence implies convergence}
Let $\Gamma$ be a group which does not have an abelian subgroup of finite index.
Let $H \subset \Gamma$ be any subgroup of finite index.
Suppose that $\{\rho_k\}$ is a sequence of discrete faithful representations
of $\Gamma$ in the group of orientation preserving isometries of $\hy^3$
such that $\{\rho_k\}$ restricted to $H$ converges.  Then $\{\rho_k\}$
converges on the entire group $\Gamma$.
\end{lemma}
\begin{proof}  One may as well assume that $H$ is a normal subgroup of $\Gamma$,
since in any case $H$ contains a normal subgroup with finite index.

Any non-elementary group, such as $H$, contains at least three non-commuting
elements of infinite order.  [In fact, $H$ must contain a free group
on $2$ generators.]  These three hyperbolic elements determine three
distinct points, there attracting fixed points,
on $\sinfty$.  Any orientation preserving isometry $g$
is determined by its action on these three points; but if $g \in \Gamma$,
this action is determined by the action of $g$ by conjugation on $H$.

If $\{\rho_k\}$ restricted to $H$ converges, then also the attracting fixed
points of its hyperbolic elements converge.  This implies that $\{\rho_k\}$
converges on all of $\Gamma$.
\end{proof}
\bigskip

\section {The geometry of pleated surfaces}

The proof of the easy case of theorem \ref{AH(acylindrical) is compact}
given in the preceding section is insufficient in the general case
for lack of control over the behavior at $\boundary M$.  For this
purpose, we will use pleated surfaces.

A {\it pleated surface} in a hyperbolic $3$-manifold $N$ is a complete
hyperbolic surface $S$, together with a continuous map $f: S \arrow N$ which
satisfies
\begin{description}
\item[(a)]  f is isometric, in the sense that every geodesic segment in $S$
is taken to a rectifiable arc in $N$ which has the same length, and\\
\item[(b)]  for each point $x \in S$, there is at least one open geodesic segment
$l_x$ through $x$ which is mapped to a geodesic segment in $N$.\end{description}

Generally, a pleated surface is not smooth, and generally $f$ is not even
locally an embedding.  These qualities are sacrificed for the sake of
properties (a) and (b) above, which turn out to generate a useful theory.
This theory has been extensively analyzed in \refto{Th1}.  We will review
a few basic facts.

\begin{proposition}[Folded Or Flat]\label{folded or flat}
If $x \in S$ is any point in the domain of a pleated surface, then
either there is a unique geodesic segment $l_x$ through $x$ which is mapped
to a geodesic segment, or $x$ has a neighborhood which is mapped isometrically
to a portion of a hyperbolic plane in the image.
\end{proposition}
\begin{proof} Suppose $l$ and $m$ are two distinct line segments through
$x$ which are mapped to line segments.  Consider a quadrilateral surrounding
$x$, with corners on $l$ and $m$.  Since its sides cannot increase in
length, none of the angles formed by $l$ and $m$ at $x$ can increase
under $f$.  This forces them to be the same in the range
as in the domain, so therefore the four
sides of the quadrilateral are mapped to geodesics.  This in turn
forces each geodesic segment from $x$ to a point on the quadrilateral
to map to a geodesic segment, so the image is a portion of a plane.
\end{proof}

The set of points for which the germ of the
segment $l_x$ is unique is called the
{\it pleating locus} $\Pi(f)$ for the pleated surface $f$.

A {\it geodesic lamination} $\lambda$ on a hyperbolic surface $S$ is a
closed set $L$, together with a foliation of $L$ by geodesics of $S$.
One example is a simple closed geodesic; other examples can be obtained
by taking for $L$ a Hausdorff limit of
a sequence of longer and longer simple closed geodesics.

When we are not being stuffy, we will ignore the distinction
between $\lambda$ and its support $L$.  Indeed, it is an easy consequence of
proposition \ref{laminations have measure 0} below that $L$
determines $\lambda$.

\begin{proposition}[Pleating Locus Is Lamination]\label{pleating locus is lamination}
The pleating locus for a pleated surface has the structure of a geodesic
lamination
\end{proposition}
\begin{proof}
It is obvious that the pleating locus is closed, since its complement
is open.

If $y \in \Pi(f)$, and if there is any geodesic segment $m$ with one
endpoint at $y$ in a direction not on $l_y$, then we claim there is
a portion of a half-plane bounded by $l_y$ containing $m$ which is
mapped by $f$ to a portion of a half-plane.  The proof is half the
proof of \ref{folded or flat}.

With these two observations, it follows that for any $x \in \Pi(f)$,
the line segment $l_x$ extends to a complete geodesic which remains
always in the pleating locus and defines $l_y$ for any point on it.
These geodesics must be coherent, defining a foliation of $\Pi(f)$,
in virtue of the fact that they cannot cross.
\end{proof}

\begin{proposition}[Laminations Have Measure 0]\label{laminations have measure 0}
The support $L$ of a geodesic lamination $\lambda$ on any surface
of finite area has measure 0.
\end{proposition}
\begin{proof} This is due to Nielsen, who studied geodesic laminations
extensively in a somewhat disguised form.
The method is to show
that the complement of $L$ has full measure. Each component of the complement
is a hyperbolic surface, with boundary consisting of geodesics, possibly with
corners (half-cusps) like the corners of an ideal triangle and cusps.  The area
of such a region can be easily computed by the Gauss-Bonnet theorem, in terms
of its Euler characteristic and number of cusps of each kind.

On the other hand, a line field can be defined in a neighborhood of
$L$ which agrees with the tangents to the lamination.  The line field
can be extended to all of $S$ but with isolated singularities.
An index is defined for such singularities, just as for vector fields;
for a line field, the index takes half-integral values.  The total
index is the Euler characteristic of $S$; when this index is allocated
among the regions in the complement of $L$, one finds that their
total area matches the area of $S$.
\end{proof}

Later we will describe how the an ideal simplicial map on the boundary of $M$
can be made to describe a pleated surface in $N_k$, much as in the example
illustrated in figure \ref{2-dimensional ideal simplicial maps}. The proof of
\ref{AH(acylindrical) is compact} will depend on knowing that these bounding
pleated surfaces remain definitely separated: they cannot collapse like a
punctured innertube.

Mathematically, the information will be conveyed in terms of the
mapping $\proj f$ of the pleating lamination $\Pi(f)$ on the domain surface $S$
into the tangent line bundle $\proj (N)$ of the target hyperbolic manifold.
Even though $f$ itself is not usually an embedding, we will show that
under mild assumptions $\proj f$ is an embedding, and from this we
will deduce that laminations in $S$ embed uniformly in $\proj f$.

A geodesic lamination $\lambda$ on a hyperbolic surface is {\it minimal}
if the closure of any leaf is all of $\lambda$, or in other words, if
there are no proper sublaminations.  A lamination is {\it recurrent}
if the closure of any half-infinite segment of any leaf contains the leaf.
It is immediate that a minimal lamination is recurrent.  There is
a partial converse:

\begin{proposition}[Finite Minimal Decomposition]\label{finite minimal decomposition}
A geodesic lamination $\lambda$ on a complete hyperbolic surface of finite
area has only a finite number of connected components, and it is
recurrent iff each component is minimal
\end{proposition}
\begin{proof}  The complement of $\lambda$ is a hyperbolic surface with
geodesic boundary, sometimes with half-cusps along its boundary.
The number of boundary components of complementary pieces is bounded
in terms of the area of $S$, so the number of components of $\lambda$
is bounded.

By elementary topology, each component of $\lambda$ has at least one
minimal set.  Suppose that $\mu$ is any minimal set of $\lambda$.
There is some $\epsilon$ such that any geodesic of $\lambda$
which comes within $\epsilon$ of the support of $\mu$ must tend toward
$\mu$ in one direction or the other:  simple geodesics within $\epsilon$
of a closed boundary component of $S - \mu$ are trapped into
a spiral around $\mu$, and simple geodesics close to a boundary component
of $S-\mu$ which has half-cusps are trapped into some half-cusp.
Therefore, if the component of $\lambda$ containing $\mu$ is
recurrent, it consists only of $\mu$.
\end{proof}

A geodesic lamination on a surface may happen to wander around only
on a small part of the surface.  For any lamination $\lambda$ on
$S$, denote by $S(\lambda)$ the smallest subsurface with geodesic boundary
containing $\lambda$.  As a special case, if $\lambda$ is a simple closed
curve, define $S(\lambda)$ to be the curve itself.  If the components of
$\lambda$ are $\mu_1$, $\mu_2$, \dots ,$\mu_n$, then $S(\lambda) = \cup S(\mu_i)$,
and the $S(\mu_i)$ have disjoint interiors.

\begin{theorem}[Laminations Cover]\label{laminations cover}
Let $f: S \arrow N$ be a map of a hyperbolic surface of finite area,
not necessarily connected, to a hyperbolic manifold, which takes any cusps
of $S$ to cusps of $N$, and which is injective on fundamental group.

Suppose that $\lambda$ is a recurrent geodesic lamination on $S$ such that
$f$ takes each leaf of $\lambda$ to a geodesic of $N$ by a local
homeomorphism (as with the pleating locus for a pleated surface).

Let $g: \lambda \arrow \proj(N)$ be the canonical lifting, and let
$\mu_1$, \dots, $\mu_n$ be the components of $\lambda$.  Then either
\begin{description}
\item[(a)] $g$ is an embedding, or \\
\item[(b)] The map $f$ restricted to the disjoint union of the surfaces
$S(\mu_i)$ factors up to homotopy through a covering map
$$\rho: \cup_. S(\mu_i) \arrow R$$
and $\cup \mu_i$ is the pullback by $p$ of a lamination $\mu$ on
$R$ which does embed in $\proj ( N )$.
\end{description}
\end{theorem}

Case (a) is actually contained in case (b); it is stated separately because it
is the typical case.  In particular, note that when $S$ is the boundary of an
acylindrical manifold, (b) cannot occur for topological reasons.  Our current
application will be in the case of a pleated surface, but the surface $S$ is
really only present to control the topology of the map of $\lambda$ into $N$.

In many ways, minimal geodesic laminations behave like simple closed geodesics;
this theorem is an illustration. In the case that $\lambda$ is a union of
simple closed geodesics, this theorem is quite clear.  

\begin{proof} 
Consider any two points $x_1, x_2 \in \lambda$ that have the same image under
$g$. Then there are geodesics $l_i: \reals \arrow S$ of $\lambda$ through $x_i$
such that $g \circ l_1 = g \circ l_2$.

Consider the diagonal lamination $\lambda \times \lambda \subset S \times S$.
The product map $l_1 \times l_2: \reals \arrow S \times S$ is a leaf of this
lamination. Since $\lambda \times \lambda$ is compact, there is a neighborhood
$U \subset S \times S$ of arbitrarily small diameter which it enters infinitely
often.  If $t$ and $s$ are any two times which $l_1 \times l_2$ enters $U$,
closed loops can be constructed, based at $x_1$ and at $x_2$: $$ \alpha_i =
l_i[0,t] * \epsilon_i * l_i \inverse [0,s] \qquad i = 1, 2,$$ where
$\epsilon_i$ is a short arc connecting $l_i(t)$ to $l_i(s)$.

The loops $g \circ \alpha_i$ are based at the point $f(x) = f(y)$ in $N$, and
they are nearly identical.  If the diameter of $U$ is less than half the
minimum value for the injectivity radius of $N$ in the image of $\lambda$, then
these loops are necessarily homotopic in $N$.

We will now prove that $g$ is a local embedding of $\lambda$ in $\proj(N)$.
Suppose that $x$ and $y$ are close to each other, and let $\delta$ be an arc on
$S$ connecting them. Then the loops $\delta * \alpha_2$ and $\alpha_1$ are
homotopic when mapped to $N$, therefore they must be homotopic in $S$.  This
implies that $l_1[0,t]$ is close to $l_2[0,t]$ for its entire length.

This works in both the forward and backward directions, varying the choice of
the points $s$ and $t$, for arbitrarily long segments of $l_x$ and $l_y$. 
Consequently, the two geodesics and the two points are identical. This proves
that $g$ is locally an embedding.

Next we will show that $g$ acts as a covering projection to its image. Define
an integer-valued function $n$ on $\lambda$ which gives the number of points in
$\lambda$ which have the same image as $x$ in $\proj(N)$.  Since $\lambda$ is
compact and since $x$ is not identified with any nearby point, $n(x)$ is
bounded.

The function $n$ must satisfy the semicontinuity condition $$ n(\lim_i(x_i))
\le \lim_i(n(x_i))$$ because points identified with $x_i$ remain well spaced
out, and any limit point of the sets identified with $x_i$ is also identified
with $x_i$. Therefore, the set $\{x | n(x) \ge m$ is closed.  It is also
saturated by $\lambda$, so it this set is a union of components.  It follows
that $n$ is constant on each $\mu_i$.

It is now easy to show that $g$ is a covering projection; this is similar to
the proof that a continuous 1-1 map between compact spaces is a homeomorphism.

What remains to be proven is that the covering projection of $\mu$ extends to
more of the surface, in fact, to the disjoint union of the $S(\mu_i)$.

The assertion is trivial in the case of any of the $\mu_i$ which are closed
geodesics, so we may assume that $\lambda$ has no such components.

In the remaining lamination, there is a finite set of special leaves of
$\lambda$, those leaves where $\lambda$ accumulates from only one side.  These
leaves can be identified topologically from $\lambda$ alone, becaue each end of
such a leaf is asymptotic to another end of such a leaf; together they form the
boundary of a half-cusp of one of the complementary regions. There are no other
pairs of leaf-ends of $\lambda$ which are asymptotic. From this, it follows
that identifications in the covering projection of  $\lambda$ extend to a
neighborhood of $\lambda$ in $S$, so that $f$ factors through a covering
projection at least on a neighborhood of $\lambda$.

Consider now any component $A$ of $S_{\mu_i} - \mu_i$. One boundary component
of $A$ is made of a chain of geodesics of $\mu_i$, the forward end of each
asymptotic to the backward end of the next, which maps as a $k$-fold covering
to its image, for some $k$.   There are two possibilities for $A$: either it is
an annulus, with another boundary component a boundary component of $S_{\mu_i}$
(or a cusp), or it is a disk.

If $A$ is an annulus, there is no problem: the covering projection defined one
of its boundary components automatically extends.

If $A$ is a disk, then the map factors up to homotopy through a branched
covering to some disk, with one branch point at the center. (Keep in mind here
that since $N$ is hyperbolic, any element of finite order in $\fund N$ is
trivial).

This means the map from the entire component $S(\mu_i)$ factors through
a branched covering.  But if the branched covering had actual
branch points, the map could not be injective on $\fund S$, contrary
to hypothesis.
\end{proof}

Often there are topological conditions that rule out the possibility
of the map of a surface into a three-manifold
factoring through any non-trivial covering, even on a union of subsurfaces.

If $S$ is a hyperbolic surface of finite area and if
$f: S \arrow N$ is a map to a hyperbolic $3$-manifold which takes
cusps to cusps, then $f$ is
{\it incompressible} if
\begin{description} \item[(a)]  $\fund f$ is injective.
\end{description}
The map $f$ is {\it doubly incompressible} if in addition
\begin{description}
\item[(b)]  Homotopy classes of maps
 $(I, \boundary I) \arrow (S, \cusps(S)$
relative to cusps map injectively to homotopy classes of maps
$(I, \boundary I) \arrow (N, \cusps N)$.\\
\item[(c)]  For any cylinder $c: S^1 \times I \arrow N$ with a factorization
of its boundary $\boundary c = f \circ c_0: \boundary ( S^1 \times I) \arrow S$
through S, if $\fund c$ is injective then either the image of
$\fund {c_0}$ consists of parabolic elements of $\fund S$, or $c_0$ extends
to a map of $S^1 \times I$ into $S$.\\
\item[(d)] Each maximal abelian subgroup of $\fund S$ is mapped to a maximal
abelian subgroup of $\fund N$.
\end{description}
We will use a slightly weaker condition,
{\it weak double incompressibility}
on surfaces where (d) is replaced by
\def\dprime{($\hbox{\rm d}'$)}
\begin{description}
\item[(d')]
Each maximal cyclic subgroup of $\fund S$ is mapped to a maximal cyclic subgroup of $\fund N$.
\end{description}
Condition (d) implies in particular that no curve on $S$ is homotopic to
a $Z^2$-cusp, while \dprime merely says that the image of
$\fund S$ is indivisible in $\fund N$.

\begin{theorem}[Laminations Inject]\label{laminations inject}
If $f: S \arrow N$ is weakly doubly incompressible, and if $\lambda$ is any
geodesic lamination on $S$ such that $f$ takes each leaf of $\lambda$ to
a geodesic in $N$ by a local homeomorphism, then the canonical lifting
$g: \lambda \arrow \proj N$ is an embedding.
\end{theorem}
\begin{proof}
From conditions (c) and \dprime of the definition of weak double
incompressibility, it easily follows that $f$ restricted to a
subsurface with geodesic boundary or to a closed geodesic of $S$
can never factor through a covering projection.
By \ref{laminations cover}, it follows that $g$ is injective
at least on the recurrent part $\rho$ of $\lambda$.

Any end of a leaf of $\lambda - \rho$ must leave all compact subsets
of $S - \lambda$.  The only possibilities are that the end either tends toward
a cusp of $S$ or a half-cusp of $S - \rho$, or spirals around a closed
geodesic of $\rho$.

\def\Cinf{S_\infty^1}
Let $\tilde \rho \subset \tilde \lambda$ be the inverse images of
$\rho$ and $\lambda$ on the universal cover $\tilde S$.
The endpoints of leaves of $\tilde \lambda$ on $\Cinf$ are either
endpoints of leaves of $\tilde \rho$, or cusps.  We will show next that
$\tilde \lambda$ maps injectively, by showing that these endpoints map
injectively to the sphere at infinity for $N$.

The set of parabolic fixed points for $\fund S$ maps injectively to the set of
parabolic fixed points for $\fund N$, by the hypothesis (b) of injectivity of
the relative fundamental groups.

A parabolic fixed point cannot be identified with an endpoint of a leaf
of $\tilde \rho$ when mapped to $\tilde N$, in virtue of the fact that the
injectivity radius for $N$ is bounded below in the image of $\rho$.

Consider any two leaves $l_1$ and $l_2$ of $\tilde \rho$ whose images in $\proj
\tilde N$ are asymptotic.  Since $\rho \arrow \proj N$ is a covering map to its
image and since $\rho$ is compact, the map $\tilde \rho \arrow \rho \arrow
\proj N$ is also a covering map to its image. It lifts to a map $\tilde \rho
\arrow \proj \tilde N$, also a covering map to its image.  This final covering
must be trivial, on account of condition (c). Comparing it to the map $\rho
\arrow \proj N$, we see that there must a fixed $\epsilon$ such that every
$\epsilon$-disk in $\rho$ maps homeomorphically to its image (since $\rho$ is
compact).

This says that $l_1$ and $l_2$ are already asymptotic in $\rho$.

It follows that the map from $\tilde \lambda$ to $\proj \tilde N$ is an
embedding.  In particular, the map $\lambda \arrow \proj N$ is a local
embedding. Define an integer-valued function $n$ on $\lambda$ which says how
many points of $\lambda$ have the same image in $\proj N$.   As before, we have
$$ n(\lim_i(x_i)) \le \lim_i(n(x_i)).$$
If a leaf of $\lambda$ limits on $\rho$, then $n$ is $1$ on the leaf, since $n$
is identically $1$ on $\rho$.  The other possibility is that the leaf limits on
cusps of $S$ at both ends; $n$ takes the value $1$ in this case because of
condition (b).
\end{proof}

We will now generalize this injectivity result to a uniform injectivity
theorem, by considering geometric limits.  For this, we will assume that the
surfaces $f: S \arrow N$ are pleated surfaces. In fact, such a pleated surface
always exists for any lamination $\lambda$ which can be mapped geodesically to
$N$, but in order to keep the logic straight, we will assume that the pleated
surface is given to us.

The main case of the uniform injectivity theorem will be when $\lambda$ is the
pleating locus for the surface, but the theorem needs to be stated with greater
generality to handle the case that the surface does not actually bend at
$\lambda$.

\begin{theorem}[Uniform Injectivity]\label{uniform injectivity}
Let $\epsilon_0 > 0$ and $A>0$ be given.
Among all doubly incompressible pleated surfaces
$$ f: S \arrow N$$
of area not greater than $A$ and laminations $\lambda \subset S$ which
are mapped geodesically by $f$, the maps
$$ g: \lambda \arrow \proj(N)$$
are uniformly injective on the $\epsilon_0$-thick part of $S$.  That is,
for every $\epsilon>0$ there is a $\delta>0$ such that for any such $S$,
$N$, $\lambda$ and $f$ and for any two points $x$ and $y \in \lambda$
whose injectivity radii are greater than $\epsilon_0$, if $d(x,y) \ge \epsilon$
then $d(g(x), g(y)) \ge \delta$.
\end{theorem}
\begin{proof}
We will make use of the geometric topology on hyperbolic manifolds equipped
with a base point and an orthogonal frame at that point.
Two manifolds are close in the geometric topology if there is a diffeomorphism
which is an approximate isometry from a ball of radius at least $R$ in
one manifold to the other manifold, taking base frame to base frame.
No condition is given on the fundamental group.

Here is a basic fact about the geometric topology:
\begin{proposition}[Geometric Topology Is Compact]\label{geometric topology is compact}
For any $\epsilon_1 > 0$, the space of hyperbolic manifolds with base
frame at a point whose injectivity radius is at least $\epsilon_1$ is
compact in the geometric topology.
\end{proposition}

We can also put a geometric topology on pleated surfaces $f:S \arrow N$
in hyperbolic manifolds, where base frames are assumed to match up
under $f$.
The pleated surface $f$ is close to $f': S' \arrow N'$
in the geometric topology if there
are approximate isometries on balls of radius $R$
from $S$ to $S'$ and from $N$ to $N'$ in balls which approximately
conjugate $f$ to $f'$.

\begin{proposition}[Pleated Surfaces Compact]\label{pleated surfaces compact}
The space of incompressible pleated surfaces $f: S \arrow N$
in hyperbolic $3$-manifolds $N$, with base frame in $S$ having
injectivity radius $\ge \epsilon_1$, is compact in the
geometric topology.
\end{proposition}

\begin{note} Without the hypothesis of incompressibility, one
would need a hypothesis on the injectivity radius of $N$.
\end{note}

\begin{proof} The base point $x_0$ of $S$ has two loops through it of length not
exceeding some constant $C$ which depends on $A$ and on $\epsilon_1$, such
that the two loops generate a free group of rank $2$.

The image $f(x_0)$ therefore also has a pair of loops through it with the
same property; it follows that the injectivity radius of $N$ at $f(x_0)$
is greater than some constant $\epsilon_2 > 0$.

The sets of potential surfaces $S$ and potential $3$-manifolds $N$ are therefore
compact in the geometric topology.  Since the maps $f$ are Lipschitz with
Lipschitz constant $1$, it follows that in any infinite sequence
$f_i: S_i \arrow N_i$
of pleated surfaces there would be a subsequence where the domain surface,
the range $3$-manifold, and the maps $f_i$ would converge to an object 
$f: S \arrow N$.  The only thing that remains to check is that the limit which
is actually a pleated surface.

If $x \in S$, and if $x$ is not an accumulation point of the $\lambda_i \subset S_i$,
then $f$ in some neighborhood of $x$ is a limit isometries of a uniform-size
hyperbolic neighborhood to a hyperbolic plane in $N_i$, so $f$ is flat in
a neighborhood of $x$.  Otherwise, geodesic arcs in $S$ which are mapped
to straight arcs by $f$ can be found as limits of arcs on $\lambda_i$.

It remains to show that the arc length of a geodesic arc $\alpha \subset S$ is
preserved by $f$. We may assume that the $\lambda_i$ converge. It is clear that
the length of the image of $\alpha$ cannot be greater than the length of
$\alpha$, since the approximations have this property. On the other hand, any
geodesic are $\alpha$ can be approximated by a geodesic arc whose length is the
sum of the lengths of its intersections with the complement of $\lambda$, since
$\lambda$ has measure 0 on $S$. Each subarc in the complement of $\lambda$ is
mapped to a straight line, so certainly its arc length is preserved; therefore,
the length of the image is at least as much as the length of $\alpha$.
\end{proof}

Suppose we are given a sequence of doubly incompressible pleated surfaces
$f_i: S_i \arrow N_i$, together with geodesic laminations $\lambda_i$
mapped geodesically, and points $x_i$ and $y_i$ on $\lambda_i$ satisfying
$\inj(x_i), \inj(y_i) \ge \epsilon_0$ and $d(g_i(x_i), g_i(y_i)) \arrow 0$.
Theorem \ref{uniform injectivity} will follow when we show that
$d(x_i, y_i) \arrow 0$.

By proposition \ref{pleated surfaces compact} there is a subsequence such
that the pleated surfaces converge with base frames chosen at either
$x_i$ or $y_i$ converge, where the first vector in the frame is tangent
to the lamination.

There is a further subsequence so that the laminations $\lambda_i$ converge
in the Hausdorff topology.  Let $f: S \arrow N$ be the limit pleated surface,
$\lambda$ the limit lamination, and $g: \lambda \arrow \proj(N)$ the canonical
lifting of $f$ restricted to $\lambda$.

\begin{lemma}[Limit Doubly Incompressible Weakly]\label{limit doubly incompressible weakly}
The limit pleated surface $f: S \arrow N$ is weakly doubly incompressible.
\end{lemma}

\begin{proof}[Proof of lemma]
We have already seen that condition (a) (injectivity)
of the definition holds, by proposition
\ref{non-elementary limit is discrete}.
\def\sithin{{(S_i)}_\thin}
\def\nithin{{(N_i)}_\thin}

We will prove (b) (injectivity of the relative fundamental group)
by showing that the maps
\def\paths{{\mathcal S}}\def\pathn{{\mathcal N}}
$$\paths =
[(I, \boundary I): (S_i, \sithin )]
\arrow
[(I, \boundary I): (N_i, \nithin )]
= \pathn$$
are injective.  Let $\alpha$ and $\beta$ be two arcs on $S_i$ with ends in
$\sithin$ representing two elements of $\paths$. Suppose that they are mapped
to the same element of $\pathn$. This means that there are arcs $\gamma$ and
$\delta$ in $\nithin$ such that
$f(\alpha)*\gamma*f(\beta\inverse)*\delta\inverse$ is null-homotopic in $N$.

If $\alpha(0)$ and $\beta(0)$ are in the same component of $\sithin$,
then they can be connected by an arc $\delta'$. Thus
$f(\delta') * \delta\inverse$ is a loop in $\nithin$. 
By condition (d) (maximal abelian groups preserved) for $f_i$, it is possible
to modify $\delta'$ to make this element trivial, so that $\delta$ is homotopic
to $f(\delta')$.  On the other hand, if $\alpha(0)$ and $\beta(0)$
are in different components of $\sithin$, then these components must
be cusps, by condition (c) (no cylinders except in cusps) for $f_i$.
A similar argument applies to $\alpha(1)$ and $\beta(1)$.

There are now various cases, all easy.  If the $\alpha(0)$ is in the same
component as $\beta(0)$ and $\alpha(1)$ is in the same component as $\beta(1)$,
then one has a closed loop in $S$ to represent the ``difference'' of the
relative homotopy classes, which maps to a trivial loop in $N$. By condition
(a) for $S$, $\alpha$ and $\beta$ must represent the same class in $\paths$. 
If all endpoints are in cusps, then $\alpha$ and $\beta$ must represent the
same class in $\paths$, by (b).  Finally, if $\alpha(0)$ and $\beta(0)$ are in
cusps and $\alpha(1)$ and $\beta(1)$ are in the same component of $\sithin$
(possibly after relabeling to reverse endpoints), then the arc $\alpha *
\delta' * \beta\inverse$ maps to a trivial element of $[(I, \boundary I),
\nithin]$, so $[\alpha] =  [\beta]$ in $\paths$.

Next we check condition (c).  Suppose that we have an incompressible cylinder
$g: S^1 \times I \arrow N$ in the limit manifold, with a factorization $g_0$ of
its boundary through $f$.  Since arbitrarily large compact subsets of $N$ are
approximately isometric to subsets of $N_i$ for large $i$, we obtain similar
cylinders $C_i$ in $N_i$, and by a small homotopy we can make their boundaries
factor through $f_i$.  By condition (c) for $f_i$, either the boundary
components go to parabolic elemenents of $S_i$, in which case the same property
holds in the limit, or there is a cylinder on $S_i$ with the same boundary. 
Consider the covering spaces of $S_i$ whose fundamental group agrees with that
of the cylinder.  If these covering spaces have core geodesics with length
going to 0, then the boundary components of $C$ are parabolic, and (c) is
satisfied.  Otherwise, the core geodesics have length bounded away from zero;
two curves of bounded length representing the generator of the fundamental
group of such a cylinder are connected by a homotopy of bounded diameter.
Therefore, there is a homotopy between the limiting curves, so (c) is satisfied
in the limit.

As for condition \dprime, suppose that $\alpha$ is a non-trivial element of
$\fund S$, and that some power $k$ of $\beta \in \fund N$  gives
$f_*(\alpha)$.  We represent this geometrically by loops $a \subset S$ and $b
\subset N$, together with a cylinder $C$ giving the homotopy from $b^k$ to
$f(a)$.  This configuration can be pushed back to the approximations $S_i$,
$N_i$ and $f_i$. It follows, by condition \dprime for $f_i$, that there is a
loop $c_i$ on $S_i$ such that $c_i^k \simeq  a_i$, where $a_i$ is the
approximation to $a$ on $S_i$.  By considering the cylindrical covering of the
hyperbolic surface $S_i$ which corresponds to $c_i$, we see that $c_i$ can be
taken to have a length less than that of $a_i$ and to lie in a neighborhood of
$a_i$ which has a diameter bounded in terms of $k$ and the length of $a_i$. 
Passing to the limit, we see that $\alpha$ is also a $k$th power.  This
verifies \dprime.
\end{proof} 

To complete the proof of theorem \ref{uniform injectivity}, we see by the lemma
that we can apply theorem \ref{laminations inject} to our limit of a sequence
of pleated surfaces in which $x_i$ and $y_i$ are mapping to points in $\proj
N_i$ which are closer and closer.  We conclude that the limit points $x =
\lim_i x_i$ and $y = \lim_i y_i$ must be equal, since their images in $\proj N$
are equal.  Therefore, $d(x_i, y_i) \arrow 0$.  This completes the proof of
\ref{uniform injectivity}.
\end{proof} 

\begin{remark}
For this proof to work, we needed both the notion of double incompressibility
and of weak double incompressibility.

In fact, the geometric limit of a sequence of doubly incompressible surfaces is
often not doubly incompressible: what happens is that a $\integers$-cusp may
enlarge to a $\integers^2$ cusp, and a limit surface can well be only weakly
doubly incompressible.

Furthermore, there are hyperbolic $3$-manifolds which admit a sequence of
weakly doubly incompressible surfaces, all homotopic, such that laminations
carried by these surfaces are {\it not} uniformly injective. The proof breaks
down because $\fund {S, S_{\thin}} \arrow \fund { N, N_{\thin}}$ is not
injective.
\end{remark}

\section {Proof that AH(M) is compact}
\label{Proof that AH(M) is compact}

Now that we are equipped with some information about the geometry
of pleated surfaces, we are prepared to return to the proof of the main theorem,
\ref{AH(acylindrical) is compact} in the general case.

We suppose that are given an acylindrical manifold $M$, together with a
sequence of elements of $AH(M)$ that is, homotopy equivalences
$h_k: M \arrow N_k$ to hyperbolic manifolds.  We will extract a subsequence
which converges.

The first step is to replace $h_k$ by an ideal simplicial map
$f_k$ such that for each component of $\boundary M$, $f_k$ describes
a pleated surface in a certain sense.  Since $f_k$ cannot be defined
on all of $S$, it is not possible for $f_k$ to be a pleated surface
literally.  Instead, we will construct pleated surfaces $p_k: S_k \arrow N_k$ together with a factorization of $f_k | \boundary M$ as a composition
of an ideal simplicial map to $S_k$ with the pleated surface,
$$f_k | \boundary M = p_k \circ i_k.$$

Triangulate $M$ with a triangulation $\tau$ which has exactly one vertex on
each boundary component, and so that the edges that lie on $\boundary M$
represent distinct non-trivial homotopy classes. For each $k$ and for each
boundary component $B_i$, there is at least one edge of $\tau$ on the component
which goes to a hyperbolic element under $h_k$.  Passing to a subsequence if
necessary, we may assume that this edge $e_i$ is the same for all $k$.

We will construct the pleated surface $p_k$ by spinning, much as in figure
\ref{2-dimensional ideal simplicial maps}.  Begin with a piecewise-straight map
of $\tau$ to $N_k$ such that the vertex of $\tau$ which is on a boundary
component of $M$ maps to a point on the geodesic homotopic to the loop formed
by $e_i$. A metric is induced on $B_i$: the distance between two points is
defined to be the infimum of the length of images in $N_k$ of paths joining the
two points.  This metric is a hyperbolic metric on $B_i$ except for one cone
point, where the total angle is greater than $2\pi$, or in other words, where
negative curvature is concentrated.

Now homotope through piecewise-straight maps by pushing the vertex around the
closed geodesic of $e_i$.  A family of metrics is obtained in this way, all
hyperbolic with one cone point.  The cone angle is bounded above (in light of
the Gauss-Bonnet theorem), and in the limit it tends toward $2\pi$.  (Most
angles of triangles tend toward $0$, except in the two triangles which have
$e_i$ as an edge; an angle of each of these tends toward $\pi$).  All the
triangles except the two degenerating triangles which border $e_i$ converge in
shape to ideal triangles, since fixed points of non-commuting hyperbolic
elements in a discrete group are distinct.  The gluing maps between ideal
triangles also converge along edges which do not border one of the two
degenerating triangles. In the limit, as the degenerating triangles get long
and thin, two of their edges are glued together.

This yields an assemblage of ideal triangles glued edge to edge homeomorphic to
$B_i - e_i$.  The resulting hyperbolic surface is mapped isometrically to
$N_k$.  The hyperbolic surface so constructed is not complete, however; in the
metric completion, two boundary curves are added, whose length is the length of
the closed geodesic of $N_k$ homotopic to $e_i$. Under the map to $N_k$, the
two boundary components are glued together. Performing the same operation in
the domain, we get a hyperbolic surface homeomorphic to $B_i$ mapped as a
pleated surface $p_k$ into $N_k$. We can choose a homeomorphism from $B_i$ to
the hyperbolic surface so that the composition into $N_k$ is isotopic to $h_k$,
enabling us to think of the hyperbolic surface as a hyperbolic structure on
$B_i$.

This construction gives the promised factorization of $f_k | \boundary M = p_k
\circ i_k$ as an ideal simplicial map followed by a pleated surface. We have a
fixed lamination $\lambda$ of $\boundary M$, such that the pleating locus is
always a sublamination of $\lambda$; the hyperbolic structure on $\boundary M$
varies with $k$.

We return to the situation as described in \S\ref{Shapes and limits},
with a good submanifold $G$ and its internal boundary $\kappa$ which
encode the way in which 3-simplices of $\tau$ are degenerating in
the sequence of ideal simplicial maps $f_k$.  We also have the
projection $p: \kappa \arrow \gamma$ describing the pattern of collapsing
of $\kappa$.  Recall that $\gamma_0$ is the union of cells whose
preimage has non-trivial fundamental group.

We will proceed in a similar way to
\S\ref{The deformation space for closed manifolds is compact},
although for the final argument we will analyze the area
of the hyperbolic structures on $\boundary M$ rather than the volume
of $N_k$.

For each edge of $\gamma - \gamma_0$ whose preimage in $\kappa$ is a cylinder,
we add a $2$-handle to obtain a new manifold $G_1$. Define $\kappa_1 \subset
\boundary G_1$ be obtained from $\kappa \cup G_1$ by adding the two bounding
disks of each of these $2$-handles, and also deleting a neighborhood of the
core transverse arc for each component of $p\inverse (\gamma - \gamma_0)$ which
is a thin strip. If you prefer, you can picture adding $2$-handles and
semi-$2$-handles respectively to the homotopically trivial thin cylinders and
thin strips.

The two boundary components of a thin strip of $\kappa$ are very close to edges
of $2$-simplices (after isotopies to make the image area small), so that they
map by $i_k$ to paths which are near edges of ideal triangles of $\boundary M -
\lambda$, hence close to leaves of $\lambda$.  Therefore there are points on
$\lambda$ near the two sides of the strip which are mapped to nearby points in
$\proj N_k$. These points can be taken to be in the thick part of the
hyperbolic structure on $\boundary M$, since thick parts of ideal triangles are
automatically contained in thick parts of the hyperbolic surface.  (Remember
that since $\kappa$ is on the boundary of the good submanifold, the
displacements of triangles which this thin strip intersects converge to bounded
values, so we can choose nearby points on the two edges of the strip which both
map to the thick parts of their triangles --- provided the $\epsilon$ for
defining thick is small enough).

Since $M$ is acylindrical, $\boundary M$ is doubly incompressible.
It follows from theorem \ref{uniform injectivity} that the two edges of
each of these strips are close together in the domain of the pleated surface,
that is, in the $k$th hyperbolic structure on $\boundary M$.  In particular, the
transverse arcs to these strips are homotopic to (unique) arcs on $\boundary M$.
Thus, we can construct a map $g_1: G_1 \arrow M$ to take
$\boundary G_1 - \kappa_1$
to $\boundary M$, and so that the image $f_k(\kappa_1)$ in $N_k$ has small area.

For any component $C$ of $\kappa_1$, $g_{1 *} \fund C$ is abelian.  If $C$ is
closed, glue on a $3$-manifold $M_C$ and map it to $M$, as we did in
\S\ref{The deformation space for closed manifolds is compact}.
If $g_{1 *} \fund {\boundary C} \subset \fund {\boundary M}$
is trivial, add a $2$-handle to each component of $\boundary C$, map it to
$\boundary M$, and finish as in the closed case.

In the remaining case, there is some component of $\boundary C$ which is
non-trivial in $\fund M$.  By acylindricity, $g_{1*}\fund C$ can only
be $\integers$, and $C$ is homotopic {\it rel} $\boundary C$ to $\boundary M$.
Let $M_C = C \times I$, and map it into $M$ to realize such a homotopy.
Define $G_2 = G_1 \cup_C M_C$, and let
$$g_2: (G_2, \boundary G_2) \arrow (M, \boundary M)$$
be a map extending $g_1$. As in \S\ref{The deformation space for closed
manifolds is compact}, $G_2$ has the property that the sequence of
representations of its fundamental group converge.

We will prove that the degree of $g_2$ is $1$ by proving the equivalent
statement that the degree of $\boundary g_2: \boundary G_2 \arrow \boundary M$
is $1$.  This degree is the same for each component of $\boundary M$;

\begin{figure}[htb]
\centering
\includegraphics{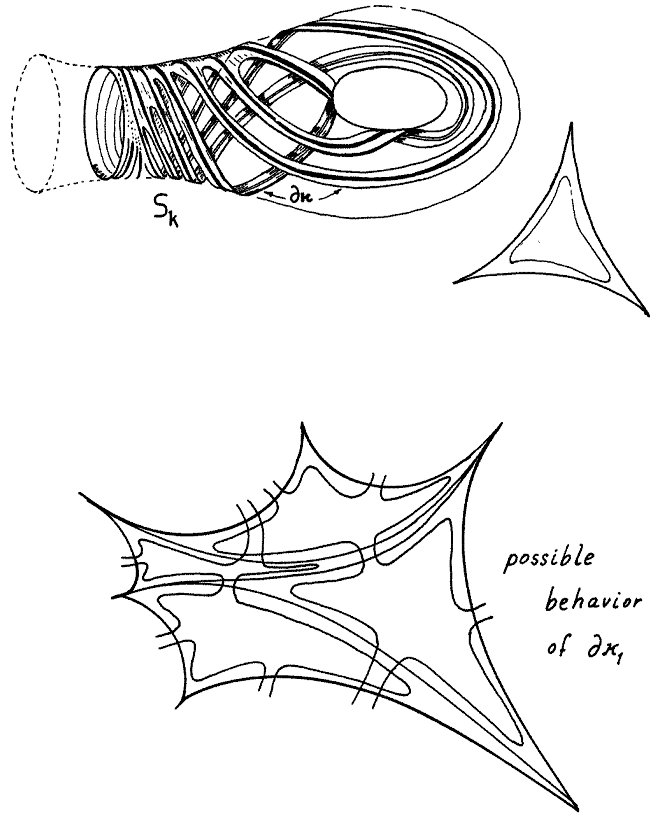}
\caption{Boundary $\kappa$. Top. The boundary of $\kappa$ has one component
going around the boundary of each
non-degenerate ideal triangle of $S_k$.  Bottom. The boundary of 
$\kappa_1$ is obtained from that of $\kappa$ by surgeries which correspond to
the thin strips of $\kappa$.  The surgeries are ``short'' because of the
uniform injectivity theorem.}
\label{boundary kappa}
\end{figure}

The original surface $\kappa \subset M$ has one boundary component for each
nondegenerate triangle of $\tau | \boundary M$.  The image by $i_k$ of each
boundary component (in the hyperbolic structure for $\boundary M$) is a curve
which goes around near the boundary of an ideal triangle of $\boundary M -
\lambda$.  The total geodesic curvature of such a curve is approximately $3
\pi$, concentrated at three left U-turns.

When we modified $G$ and $\kappa$ to obtain $G_1$ and $\kappa_1$, and
constructed the map $g_1: G_1 \arrow M$, the curves on $\boundary M$ changed.
The $2$-handles which were added did not affect the boundary curves, but when
the neighborhood of a transverse arc from each thin strip of $\kappa$ was
deleted and pushed out to the boundary (adding the ``semi-$2$-handles''), this
operation did affect the boundary curves.  The effect was to pair each long
straight segment of a boundary curve with an opposite long straight segment and
perform a 0-surgery, thereby creating two new U-turns.  The curve system after
these operations is likely not to be embedded.  The new U-turns are in an
opposite sense from the U-turns coming from corners of triangles, so that the
components of $\boundary \kappa_1$ have total geodesic curvature near 0.

We claim that each component of $\boundary \kappa$ is homotopic on $S_k$ to its
geodesic on $S_k$ by a homotopy of small mass (as $k \arrow \infty$.)

\begin{figure}[htb]
\centering
\includegraphics{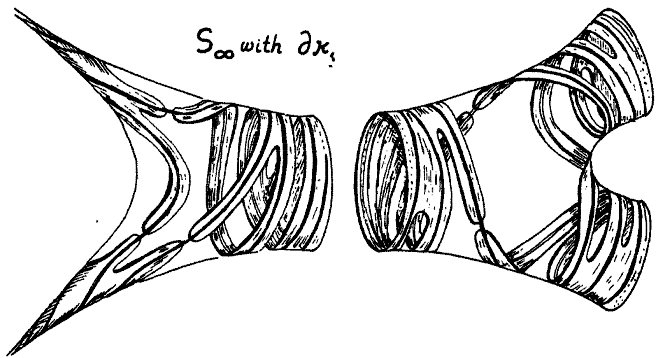}
\caption{An alternate set of identifications for the ideal triangles of $S_k$
is approximately determined by the geometry of the thin strips of $\boundary
\kappa$.  This gives a new surface $R_k$ which is the interior of a surface
with geodesic boundary components; on the new surface, the system of boundary
curves of $\kappa_1$ becomes embedded.}
\label{alternate identifications}
\end{figure}

Perhaps the clearest way to establish this claim is to make use of the thin
strips of $\kappa$ to define an alternate set of identifications of the ideal
triangles of $S_k$. That is, the two sides of each of these thin strips runs
along near two sides of ideal triangles; identify these two sides by a map as
approximately determined by the geometry of the strip, to define a new surface
$R_k$ (see figure \ref{alternate identifications}) The surface $R_k$ might not
be closed, and it might not be metrically complete.  Topologically, $R_k$ is
equivalent to the surface obtained by identifying true $2$-simplices in the
given combinatorial pattern, then removing the vertices. Its metric completion
$\bar R_k$ in general is a hyperbolic surface with either a cusp or a geodesic
boundary component for each missing vertex (see \refin{pp. 3.17-3.22}{Th1}.)
Since we are free to make slight adjustments in the gluing patterns, we may
assume that all missing vertices correspond to geodesics. In $R_k$, the
components of $\boundary \kappa_1$ map to simple, disjoint curves. 
Topologically, each of these curves circles around a neighborhood of one of the
missing vertices.  By the Gauss-Bonnet theorem, the punctured disk which gives
the homology of any of these curves to its boundary geodesic has small area on
$\bar R_k$.

For large $k$, there is a map of $\bar R_k$ to $N_k$, constructed by mapping
neighborhoods of the edges of the ideal triangulation over the thin strips of
$\kappa$. This map can be made an approximate isometry without difficulty
except perhaps on the union of the thin parts of the points of ideal triangles.
By theorem \ref{uniform injectivity}, the edges of the thin strips are close to
each other on the pleated surface $S_k$, so this defines a map of $\bar R_k$ to
$S_k$ as well, also an approximate isometry except possibly in the union of the
thin parts of points of ideal triangles. Each component of the union of the
thin parts of the points of the ideal triangles is topologically an open
cylinder, and the only invariant up to approximate isometry is the translation
length of its holonomy.  Therefore, by easy adjustments, the map can be made an
approximate isometry everywhere.

Within the union of the thick parts of the triangles of $R_k$, each of the
images of $\boundary \kappa_1$ has a homology of small mass to a curve with
uniformly small curvature.  This homology pushes forward to a similar homology
on $S_k$, from $\boundary \kappa_1$ to some curve of small absolute curvature
on $S_k$. 

Any curve of uniformly small curvature on a complete hyperbolic surface of
finite area is homotopic to its geodesic or cusp by a homology whose mass is a
constant near $0$ times the length of the curve.  Putting these two homologies
together, we establish the claim.

Finally, we added certain submanifolds $M_C$ to $G_1$, yielding $G_2$, and we
extended the map $g_1$ over the $M_C$'s to construct a map  $g_2: G_2 \arrow
M$. The image of the fundamental group of each $M_C$ in $N$ is abelian, so the
map of $\boundary G_2 - \boundary g_1 =\boundary M_C - \kappa$ to $S_k$ is
homotopic to a geodesic of $S_k$. The image $2$-chain is therefore a sum of
chains as constructed above.

Consequently, the local degree of the map of $\boundary G_2$ to $S_k$
is approximately $1$, therefore exactly $1$.  It follows that the
degree of $g_2$ is $1$.

The proof of theorem
\ref{AH(acylindrical) is compact}
is completed by an application of lemma
\ref{finite index convergence implies convergence}, as
in \S\ref{The deformation space for closed manifolds is compact}.
\end{proof}

\section{Manifolds with designated parabolic loci}

Often in the study of hyperbolic manifolds, it is important to specify
information about cusps.  In this section, we will generalize theorem
\ref{AH(acylindrical) is compact} to take into account such specifications.

Data concerning cusps for hyperbolic structures on a manifold $M^n$
can conveniently be encoded by giving a submanifold
$P^{n-1} \subset \boundary M$.  $P$ should have the following properties:
\begin{description}
\item[(a)] the fundamental group of each of its components injects into the
fundamental group of $M$.\\
\item[(b)] the fundamental group of each of its components
contains an abelian subgroup with finite index.\\
\item[(c)]  any cylinder 
$$C: (S^1 \times I, \boundary S^1 \times I) \arrow (M, P)$$
such that $\fund C$ is injective is homotopic {\it rel} boundary to $P$\\
\item[(d)] $P$ contains every component of $\boundary M$ which has an abelian subgroup
of finite index.
\end{description}
A manifold $M$ equipped with a submanifold $P$ of its boundary satisfying
these conditions is a {\it pared manifold}.
Condition (d) is in a sense unnecessary, since a component of the
boundary of $M$ which has an abelian subgroup of finite index
is necessilarily parabolic for any element of $H(M)$.  For the same reason,
it is convenient to include it.
In the three-dimensional oriented case, $P$ is a union of tori and annuli.
Define $H(M,P)$ to be the set of complete hyperbolic manifolds $M$
together with a homotopy equivalence of $(M, P)$ to $(N, \cusps)$, where
we represent cusps by disjoint horoball neighborhoods.  Since
$H(M,P) \subset H(M)$, there are induced topologies $AH(M,P)$, $GH(M,P)$
and $QH(M,P)$.

A pared 3-manifold $(M, P)$ is {\it acylindrical} if
$\boundary M - P$ is incompressible and if every cylinder
$$C: (S^1 \times I, \boundary S^1 \times I) \arrow (M, \boundary M - P)$$
such that $\fund C$ is injective is homotopic {\it rel} boundary
to $\boundary M$.  With this definition, it follows that any
hyperbolic structure $N \in H(N,P)$ and for any pleated surface $f$
representing $\boundary M - P$ is doubly incompressible.  Our main
theorem generalizes to

\begin{theorem}[AH(Pared Acylindrical) Is Compact]\label{AH(pared acylindrical) is compact}
If $(M^3, P^2)$ is an acylindrical pared manifold, then AH(M,P) is compact
\end{theorem}

This proof is practically identical to the proof in
\S\ref{Proof that AH(M) is compact}, but we will repeat it anyway.
This theorem is definitely not contained in the previous version.
For example, an important case is when $M$ is a handlebody or a surface $\times$
interval.  In either case, it is easy to see that $AH(M)$ is noncompact,
yet even for these manifolds, there are many choices of
$P \subset \boundary M$ such that $(M, P)$ is an acylindrical pared manifold.

\begin{proof}
Let $\{ N_k \}$ be a sequence of elements of $AH(M, P)$. As in \ref{Proof that
AH(M) is compact}, we can choose a triangulation and ideal simplicial maps
$f_k: M \arrow N_k$ such that $f_k | \boundary M - P$ factors as an ideal
simpliecial map composed with a pleated surface, $$ f_k | \boundary M - P = p_k
\circ i_k.$$ The domain of the pleated surface $P_k$ is a complete hyperbolic
structure of finite area on $\boundary M - P$, and after passing to a
subsequence if necessary, the pleating locus is always contained in a fixed
lamination $\lambda$ which cuts the surface into ideal triangles.

Pass to a subsequence, as described in \S\ref{Shapes and limits}, so that edge
invariants for nondegenerate simplices all converge (perhaps to degenerate
values).  The good submanifold $G$ then intersects $\boundary M$ as in
\S\ref{Proof that AH(M) is compact} only in the complement of $P$.  The image
under $i_k$ of this intersection is the union of the main parts of the ideal
triangles of $\boundary M - P - \lambda$. As before, a limiting gluing diagram
for these ideal triangles is obtained by gluing according to the thin strips of
the surface $\kappa$, yielding a hyperbolic surface $S_\infty$.  By the uniform
injectivity theorem, $S_\infty$ has an approximate isometry to the $k$th
hyperbolic structure on $\boundary M - P$, for large $k$.

A manifold $G_1$ with a distinguished submanifold
$\kappa_1 \subset \boundary G_1$
is constructed just as in \S\ref{Proof that AH(M) is compact}
by gluing on $2$-handles and semi-$2$-handles to each portion of $\kappa$
which is the preimage of an edge of $\gamma - \gamma_0$.  We can map $G_1$ to
$M$ so that $\boundary G_1 - \kappa_1$ maps to $\boundary M$, and so that
only a small change is made in the area of $i_k(\boundary G_1 - \kappa_1)$.  

Each component of $\boundary \kappa_1$ now corresponds to a curve
circumnavigating a missing vertex of $S_\infty$, and the image of the
fundamental group of each component of $\kappa_1$ in the fundamental group of
$M$ is  abelian.  If a component $C$ of $\kappa_1$ is closed, it can be capped
off with a manifold $M_C$ as before.  Otherwise, $C$ is homotopic to an annulus
on $\boundary M$. A product manifold $M_C$ realizing a homotopy can be glued
on, giving a new manifold $G_2$ on whose fundamental group the sequence of
representations converges, with a map $(G_2, \boundary G_2) \arrow (M,
\boundary M)$. The only difference in the picture is that $\boundary G_2$ may
hit annuli in $P$.  Nonetheless, the area argument works precisely as before,
to show that the map $(G_2, \boundary G_2) \arrow (M, \boundary M)$ has degree
1. Therefore, $\fund {G_2} \arrow \fund M$ is surjective, so the sequence of
representations of $\fund M$ converges.  This proves that $AH(M,P)$ is compact.

\end{proof}

\bibliographystyle{alpha}
\bibliography{../hype}
\end{document}